\documentclass[12pt]{amsart}
\usepackage{amsmath,amsthm,latexsym,amscd,amsbsy,amssymb,amsfonts,fleqn,leqno,
euscript, pb-diagram}

\numberwithin{equation}{section}

\newcommand{\LC[1]}{\mathrm{LC}^{#1}}

\newcommand{\II}{\mathbb I}

\newcommand{\Ddim}{\dim_\triangle}

\newcommand{\U}{\mathcal U}
\newcommand{\V}{\mathcal V}
\newcommand{\W}{\mathcal W}
\newcommand{\St}{St}

\newcommand{\I}{\mathbb I}

\newcommand{\e}{\varepsilon}
\newcommand{\diam}{\mathrm{diam}}
\newcommand{\mesh}{\mathrm{mesh}}

\newtheorem{thm}{Theorem}[section]
\newtheorem{pro}[thm]{Proposition}
\newtheorem{lem}[thm]{Lemma}
\newtheorem{cor}[thm]{Corollary}

\newtheorem{question}{Question}

\begin{document}


\title[Approximation by light maps and parametric Lelek maps]
{Approximation by light maps and parametric Lelek maps}

\author{Taras Banakh}
\address{Department of Mechanics and Mathematics,
Ivan Franko Lviv National University (Ukraine) and \newline Instytut Matematyki,
Akademia \'Swi\c etokrzyska w Kielcach (Poland)} \email{tbanakh@franko.lviv.ua}

\author{Vesko  Valov}
\address{Department of Computer Science and Mathematics, Nipissing University,
100 College Drive, P.O. Box 5002, North Bay, ON, P1B 8L7, Canada}
\email{veskov@nipissingu.ca}
\thanks{The second author was partially supported by NSERC Grant 261914-03.}

\keywords{dimension, $n$-dimensional maps, $n$-dimensional Lelek
maps, dendrites, Cantor $n$-manifolds, general position properties}

\subjclass{Primary 54F45; Secondary 55M10}


\begin{abstract}
The class of metrizable spaces $M$ with the following approximation
property is introduced and investigated: $M\in AP(n,0)$ if for every
$\e>0$ and a map $g\colon\I^n\to M$ there exists a $0$-dimensional
map $g'\colon\I^n\to M$ which is $\e$-homotopic to $g$. It is shown
that this class has very nice properties. For example, if $M_i\in
AP(n_i,0)$, $i=1,2$, then $M_1\times M_2\in AP(n_1+n_2,0)$.
Moreover, $M\in AP(n,0)$ if and only if each point of $M$ has a
local base of neighborhoods $U$ with $U\in AP(n,0)$. Using the
properties of $AP(n,0)$-spaces, we generalize some results of Levin
and Kato-Matsuhashi concerning the existence of residual sets of
$n$-dimensional Lelek maps.
\end{abstract}

\maketitle

\markboth{}{Approximation by light maps}



\section{Introduction}
All spaces in the paper are assumed to be metrizable and all maps
continuous. By $C(X,M)$ we denote all maps from $X$ into $M$.
Unless stated otherwise, all function spaces are endowed with the
source limitation topology.

One of the important properties of the $n$-dimensional cube $\I^n$
or the Euclidean space $\mathbb{R}^n$, widely exploited in dimension
theory, is that any map from an $n$-dimensional compactum $X$ into
$\I^n$ (resp., $\mathbb{R}^n$) can be approximated by maps with
$0$-dimensional fibers. The aim of this article is to introduce and
investigate the class of spaces having that property. More
precisely, we say that a space $M$ has the {\em
$(n,0)$-approximation property} (br., $M\in AP(n,0))$ if for every
$\e>0$ and a map $g\colon\I^n\to M$ there exists a $0$-dimensional
map $g'\colon\I^n\to M$ which is $\e$-homotopic to $g$. Here, $g'$
is $\e$-homotopic to $g$ means that there is an $\e$-small homotopy
$h\colon\I^n\times\I\to M$ connecting $g$ and $g'$. It is easily
seen that this definition doesn't depend on the metric generating
the topology of $M$. If $M$ is $\LC[n]$, then $M\in AP(n,0))$ if and
only if for every $\e>0$ and $g\in C(\I^n,M)$ there exists a
$0$-dimensional map $g'$ which is $\e$-close to $g$. We show that
the class of $AP(n,0)$-spaces is quite large and it has very nice
properties. For example, if $M_i\in AP(n_i,0)$, $i=1,2$, then
$M_1\times M_2\in AP(n_1+n_2,0)$. We also prove that $M\in AP(n,0)$
if and only if each point of $M$ has a local base of neighborhoods
$U$ with $U\in AP(n,0)$. Moreover, every path-connected compactum
$M\in AP(n,0)$ is a $(V^n)$-continuum in the sense of P. Alexandroff
\cite{ps}. In particular, according to \cite{hv}, any such $M$ is a
Cantor $n$-manifold, as well as a strong Cantor $n$-manifold in the
since of Had\v{z}iivanov \cite{nh}. All complete $\LC[0]$-spaces
without isolated points are $AP(1,0)$, as well as, every manifold
modeled on the $n$-dimensional Menger cube or the $n$-dimensional
N\"{o}beling space has the $AP(n,0)$-property.

The class of $AP(n,0)$-spaces is very natural for obtaining results
about approximation by dimensionally restricted maps. We present
such a result about $n$-dimensional Lelek maps. Recall that a map
$g\colon X\to M$ between compact spaces is said to be an {\em
$n$-dimensional Lelek map} \cite{l} if the union of all non-trivial
continua contained in the fibers of $g$ is of dimension $\leq n$.
For convenience, if $n\leq 0$, by an $n$-dimensional Lelek map we
simply mean a $0$-dimensional map. Every $n$-dimensional Lelek map
$g$ between compacta is at most $n$-dimensional because the
components of each fiber of $g$ is at most $n$-dimensional (we say
that $g$ is $n$-dimensional if all fibers of $g$ are at most
$n$-dimensional). Lelek \cite{le} constructed such a map from
$\I^{n+1}$ onto a dendrite. Levin \cite{l} improved this result by
showing that the space $C(X,\I^n)$ of all maps from $X$ into $\I^n$
contains a dense $G_{\delta}$-subset consisting of
$(m-n)$-dimensional Lelek maps for any compactum $X$ with $\dim
X=m\leq n$. Recently Kato and Matsuhashi \cite{km} established a
version of the Levin result with $\I^n$ replaced by more general
class of spaces. This class consists of complete metric $ANR$-spaces
having a piecewise embedding dimension $\geq n$.

In the present paper we generalize the result of Kato-Matsuhashi in
two directions. First, we establish a parametric version of their
theorem and second, we show that this parametric version holds for
$AP(n,0)$-spaces (see Section 6 where it is shown that the class of
$AP(n,0)$-spaces contains properly the class of $ANR$'s with
piecewise embedding dimensions $\geq n$). Here is our result about
parametric Lelek maps (a more general version is established in
Section 3; moreover we also present a version of Theorem 1.1 when
$M$ has a property weaker than $AP(n,0)$, see Theorem 4.5).

\begin{thm}\label{lelek}
Let $f\colon X\to Y$ be a perfect map with $\Ddim(f)\leq m$, where
$X$ and $Y$ are metric spaces. If $M\in AP(n,0)$ is completely
metrizable, then there exists a $G_\delta$-set $\mathcal{H}\subset
C(X,M)$ such that every simplicially factorizable map in $C(X,M)$
is homotopically approximated by maps from $\mathcal{H}$  and
for every $g\in\mathcal{H}$ and $y\in Y$ the restriction
$g|f^{-1}(y)$ is an $(m-n)$-dimensional Lelek map.
\end{thm}

Let us explain the notions in Theorem 1.1. A map $g\in C(X,M)$ is
homotopically approximated by
maps from $\mathcal{H}$ means that for every function $\e\in
C(X,(0,1])$ there exists $g'\in\mathcal{H}$ which is $\e$-homotopic
to $g$. Here, the maps $g$ and $g'$ are said to be $\e$-homotopic,
if there is a homotopy $h\colon X\times\I\to M$ connecting $g$ and
$g'$ such that each set $h(\{x\}\times\I)$ has a diameter $<\e(x)$,
$x\in X$.

The function space $C(X,M)$ appearing in this theorem is endowed
with the source limitation topology whose neighborhood base at a
given function $f\in C(X,M)$ consists of the sets
$$B_\varrho(f,\e)=\{g\in C(X,M):\varrho(g,f)<\e\},$$ where $\varrho$ is a
fixed compatible metric on $M$ and $\e:X\to(0,1]$ runs over
continuous positive functions on $X$. The symbol $\varrho(f,g)<\e$
means that $\varrho(f(x),g(x))<\e(x)$ for all $x\in X$. Since $X$ is
metrizable, the source limitation topology doesn't depend on the
metric $\rho$ \cite{nk} and it has the Baire property provided $M$
is completely metrizable \cite{munkers}.

We say that a map $g\colon X\to M$ is simplicially factorizable
\cite{bv} if there exists a simplicial complex $L$ and two maps
$g_1\colon X\to L$ and $g_2\colon L\to M$ such that $g=g_2\circ
g_1$. In each of the following cases the set of simplicially
factorizable maps is dense in $C(X,M)$ (see \cite[Proposition
4]{bv}): (i) $M$ is an $ANR$; (ii) $\dim X\leq k$ and $M$ is
$\LC[k-1]$; (iii) $X$ is a $C$-space and $M$ is locally
contractible.

The dimension $\Ddim(f)$ was defined in \cite{bv}: $\Ddim(f)$ of a
map $f:X\to Y$ is equal to the smallest cardinal number $\tau$ for
which there is a map $g:X\to \mathbb I^\tau$ such that the diagonal
product $f\Delta g:X\to Y\times \mathbb I^\tau$ is a $0$-dimensional
map. For any perfect map $f:X\to Y$ between metric spaces we have:
(i) $\dim(f)\leq \Ddim(f)$; (ii) $\Ddim(f)=\dim(f)$ if $Y$ is a
$C$-space, see \cite{bp:96} and \cite{tv}; (iii) $\Ddim(f)\leq\dim
(f)+1$ if the spaces $X,Y$ are compact, see \cite{lev}.

Since every metric space admitting a perfect finite-dimensional map
onto a $C$-space is also a $C$-space \cite{hy}, Theorem 1.1 implies
the following (here, $\mathcal{H}$ is homotopically dense in
$C(X,M)$ if every $g\in C(X,M)$ is homotopically approximated by
maps from $\mathcal{H}$):

\begin{cor}
Let $f:X\to Y$ be a perfect $m$-dimensional map between metric
spaces with $Y$ being a $C$-space. If $M\in A(n,0)$ is completely
metrizable,  then in each of the following cases there exists a
homotopically dense $G_{\delta}$-subset $\mathcal{H}\subset C(X,M)$ consisting of
maps $g$ such that $g|f^{-1}(y)$ is an $(m-n)$-dimensional Lelek map
for every $y\in Y$:
\begin{itemize}
\item $M$ is locally contractible;
\item $X$ is finite dimensional and $M\in\LC[k-1]$ for $k=\dim X$.
\end{itemize}
\end{cor}

The paper is organized as follows. In Section 2 some preliminary
results about $AE(n,0)$-spaces are presented. In Section 3 we
establish a generalized version of Theorem 1.1. Section 4 is devoted
to almost $AP(n,0)$-spaces. The final two sections contain some
interesting properties of $AP(n,0)$-spaces. In Section 5 we prove
that $AP(n,0)$-property has a local nature, as well as, $M_1\times
M_2\in AP(n_1+n_2,0)$ provided each $M_i\in AP(n_i,0)$ is completely
metrizable. It is also established in this section that every path
connected compactum is a $(V^n)$-continuum. Section 6 is devoted to
the interplay of $AP(n,0)$-property and the general position
properties $m-\overline{DD}^{\{n,k\}}$-properties. In particular, it
is shown that every point of a locally path-connected
$AP(n,0)$-space $X$ is a homological $Z_{n-1}$-point in $X$. Another
result from this section states that every completely metrizable
space possessing the disjoint $(n-1)$-disks property is an
$AP(n,0)$-space.

\section{Preliminary results about $AP(n,0)$-spaces}

In this section we established some preliminary results on
$AP(n,0)$-spaces which are going to be used in next sections.

Suppose $M$ is a metric space and $\varepsilon >0$. We write
$d_n(M)<\varepsilon$ if $M$ can be covered by an open family
$\gamma$ such that $\diam U<\e$ for all $U\in\gamma$ and
$ord(\gamma)\leq n+1$ (the last inequality means that at most $n+1$
elements of $\gamma$ can have a common point). It is easily seen
that if $F\subset M$ is closed and $d_n(F)<\varepsilon$, then $F$ is
covered by an open in $M$ family $\gamma$ of mesh $<\varepsilon$ and
$ord\leq n+1$. We also agree to denote by $cov(M)$ the family of all
open covers of $M$. Let us begin with the following technical lemma.

\begin{lem}
Let $Z=A\cup B$ be a compactum, where $A$ and $B$ are closed subsets
of $Z$. Suppose $C\subset Z$ is closed such that $\dim C\cap A\leq
0$ and $d_0(C\cap B)<\varepsilon$. Then $d_0(C)<\varepsilon$.
\end{lem}

\begin{proof}
Since $d_0(C\cap B)<\e$, there exists a disjoint open family
$\gamma=\{W_1,..,W_k\}$ in $B$ of $\mesh<\e$. We extend every $W_i$
to an open set $\tilde{W_i}$ in $Z$ such that $\tilde{W_i}\cap
B=W_i$ and $\tilde{\gamma}=\{\tilde{W_i}:i=1,..,k\}$ is a disjoint
family of $\mesh <\varepsilon$ (this can be done because $B$ is
closed in $Z$). Observe that
$C_1=C\backslash\bigcup_{i=1}^{i=k}\tilde{W_i}$ is a closed subset
of $C\cap A$ disjoint from $B$. Since $\dim C\cap A\leq 0$, there
exists a clopen set $C_2$ in $C\cap A$ disjoint from $B$ and
containing $C_1$. Obviously, $C_2$ is clopen in $C$. Moreover, $\dim
C_2\leq 0$. Hence, there exists a cover $\omega=\{V_j:j=1,..,m\}\in
cov(C_2)$ consisting of clopen subsets of $C_2$ with $\diam V_j<\e$
for every $j$. Then $\omega_1=\omega\cup\{\tilde{W_i}\cap
C\backslash C_2\}$ is a disjoint open cover of $C$ and
$mesh(\omega_1) <\e$. So, $d_0(C)<\e$.
\end{proof}

The proof of next lemma is extracted from \cite[Theorem
4]{bv}.

\begin{lem}
Let $G\subset C(X,M)$, where $(M,\varrho)$ is a complete metric
space. Suppose $(U(i)_{i\geq 1}$ is a sequence of open subsets of
$C(X,M)$ such that
\begin{itemize}
\item for any $g\in G$, $i\geq 1$ and any function $\eta\in
C(X,(0,1])$ there exists $g_i\in B_\varrho(g,\eta)\cap U(i)\cap G$
which is $\eta$-homotopic to $g$.
\end{itemize}
Then, for any $g\in G$ and $\e\colon X\to (0,1]$ there exists
$g'\in\bigcap_{i=1}^{\infty} U(i)$ and an $\e$-homotopy connecting
$g$ and $g'$. Moreover, $g'|A=g_0|A$ for some $g_0\in C(X,M)$ and
$A\subset X$ provided $g_i|A=g_0|A$ for all $i$.
\end{lem}

\begin{proof}
For fixed $g\in G$ and $\e\in C(X,(0,1])$ let $g_0=g$ and
$\e_0=\e/3$. We shall construct by induction a sequence $(g_i:X\to
M)_{i\geq 1}\subset G$, a sequence $(\e_i)_{i\geq 1}$ of positive
functions, and a sequence $(H_i:X\times[0,1]\to M)_{i\geq 1}$ of
$\e_{i-1}$-homotopies satisfying the conditions:
\begin{itemize}
\item $H_{i+1}(x,0)=g_{i}(x)$ and $H_{i+1}(x,1)=g_{i+1}(x)$ for every $x\in X$;
\item $g_{i+1}\in B_\varrho(g_i,\e_i)\cap U(i+1)\cap G$;
\item $\e_{i+1}\leq\e_i/2$;
\item $B_\varrho(g_{i+1},3\e_{i+1})\subset
U(i+1)$.
\end{itemize}

Assume that, for some $i$, we have already constructed maps
$g_1,\dots,g_i$, positive numbers $\e_1,\dots,\e_i$, and homotopies
$H_1,\dots,H_i$ satisfying the above conditions. According to the
hypotheses, there exists a map $g_{i+1}\in B_\varrho(g_i,\e_i)\cap
 U(i+1)\cap G$ such that $g_{i+1}$ is $\e_i$-homotopic to $g_i$. Let
$H_{i+1}:X\times[0,1]\to M$ be an $\e_i$-homotopy connecting the
maps $g_{i}$ and $g_{i+1}$. Since the set $U(i+1)$ is open in
$C(X,M)$, there is a positive function $\e_{i+1}\leq\e_i/2$ such that
$B_\varrho(g_{i+1},3\e_{i+1})\subset U(i+1)$. This completes the
inductive step.

It follows from the construction that the function sequence
$(g_i)_{i\geq 1}$ converges uniformly to some continuous function
$g':X\to M$. Obviously, $\varrho(g',g_i)\leq \sum_{j=i}^\infty\e_j\leq
2\e_i$ for every $i$. Hence, according to the choice of the
sequences $(\e_i)$ and $(g_i)$, $g'\in U(i)$ for every $i\geq 1$.
So, $g'\in\bigcap_{i=1}^{\infty} U(i)$. Moreover, the
$\e_{i-1}$-homotopies $H_i$ compose an $\e$-homotopy
$H:X\times[0,1]\to M$
$$
H(x,t)=\begin{cases} \displaystyle
H_i\big(x,2^i(t-1+\frac1{2^{i-1}})\big)&\mbox{if $\displaystyle t\in
[1-\frac1{2^{i-1}},1-\frac1{2^{i}}]$,
$i\geq 1$};\\
g'(x)&\mbox{if $t=1$}.
\end{cases}
$$
connecting $g_0=g$ and $g'$.

If follows from our construction that $g'|A=g_0|A$ if $g_i|A=g_0|A$
for all $i\geq 1$.
\end{proof}

\begin{pro}
Let $M$ be a completely metrizable space with the $AP(n,0)$-property
and $X$ be an $n$-dimensional compactum. Then $$\mathcal E=\{g\in
C(X,M): \dim g\leq 0\}$$ is a $G_{\delta}$-subset of $C(X,M)$ such
that every simplicially factorizable map is homotopically
approximated by maps from $\mathcal E$.
\end{pro}

\begin{proof}
We fix a metric on $X$ and a complete metric $\varrho$ on $M$. For
every $\e>0$ let $C(X,M;\e)$ be the set all maps $g\in C(X,M)$ with
$d_0(g^{-1}(a))<\varepsilon$ for every $a\in M$. It is easily seen
that $C(X,M;\e)$ is open in $C(X,M)$ and $\mathcal
E=\bigcap_{i=1}^{\infty}C(X,M;1/i)$. So, $\mathcal E$ is a
$G_{\delta}$-subset of $C(X,M)$.

{\em Claim. For every simplicially factorizable map $h\in C(X,M)$
and positive numbers $\eta$ and $\delta$ there exists a simplicially
factorizable map $h'\in C(X,M;\eta)$ which is $\delta$-homotopic to
$h$}.

To prove the claim, fix a simplicially factorizable map $h\in
C(X,M)$ and positive numbers $\delta$, $\eta$. Then there exists a
simplicial complex $L$ and maps $q_1\colon X\to L$, $q_2\colon L\to
M$ with $h=q_2\circ q_1$. We can assume that $L$ is finite and
$q_2(L)\subset\cup\gamma$, where $\gamma=\{B_{\varrho}(y,\delta/4):
y\in g(X)\}$. Because $L$ is an $ANR$ and $\dim X\leq n$, there is a
polyhedron $K$ with $\dim K\leq n$ and two maps $f\colon X\to K$,
$\alpha: K\to L$ such that $f$ is an $\eta$-map  and $q_1$ and
$\alpha\circ f$ are $q_2^{-1}(\gamma)$-homotopic. So, $h$ and
$q_2\circ\alpha\circ f$ are $\delta/2$-homotopic. Moreover, there is
a cover $\omega\in cov(K)$ such that $\diam f^{-1}(W)<\eta$ for
every $W\in\omega$. Let $\theta$ be a Lebesgue number of $\omega$.
It remains to find a map $h^*\colon K\to M$ which is
$\delta/2$-homotopic to $q_2\circ\alpha$ and $h^*\in C(K,M;\theta)$
(then $h'=h^*\circ f$ would be a simplicially factorizable map
$\delta$-homotopic to $h$ and $h'\in C(X,M;\eta)$).

To find such a map $h^*\colon K\to M$, let
$\{\sigma_1,..,\sigma_m\}$ be an enumeration of the simplexes of $K$
and $K_i=\bigcup_{j=1}^{j=i}\sigma_j$. We are going to construct by
induction maps $h_i\colon K\to M$, $i=0,..,m$, satisfying the
following conditions:
\begin{itemize}
\item $h_0=q_2\circ\alpha$;
\item $h_i|K_i$ belongs to $C(K_i,M;\theta)$, $1\leq i\leq m$;
\item $h_i$ and $h_{i+1}$ are $(\delta/2m)$-homotopic, $i=0,..,m-1$.
\end{itemize}

Let $\V=\{B_{\varrho}(y,\delta/4m): y\in M\}\in cov(M)$ and assume
that $h_i$ has  already been constructed. Since $h_i|K_i\in
C(K_i,M;\theta)$, every fiber $h_i^{-1}(y)\cap K_i$ of $h_i|K_i$,
$y\in h_i(K_1)$, is covered by a finite open and disjoint family
$\Omega(y)$ in $K_i$ with $\mesh(\Omega(y))<\theta$. Using that
$h_i|K_i$ is a perfect map, we find a cover $\V_i\in cov(M)$ such
that $\V_i$ is a star-refinement of $\V$ and if $\St(y,\V_i)\cap
h_i(K_i)\neq\emptyset$ for some $y\in M$, then there is $z\in
h_i(K_i)$ with $h_i^{-1}(\St(y,\V_i))\cap K_i\subset\Omega(z)$.
Since $M$ has the $AP(n,0)$-property, there exists a $0$-dimensional
map $p_i\colon\sigma_{i+1}\to M$ which is $\V_i$-homotopic to
$h_i|\sigma_{i+1}$. By the Homotopy Extension Theorem, $p_i$ can be
extended to a map $h_{i+1}\colon K\to M$ being $\V_i$-homotopic to
$h_i$. Then $h_{i+1}$ is $(\delta/2m)$-homotopic to $h_i$. To show
that $h_{i+1}|K_{i+1}$ belongs to $C(K_{i+1},M;\theta)$, we observe
that $h_{i+1}^{-1}(y)\cap K_{i+1}=\big(h_{i+1}^{-1}(y)\cap
K_i\big)\cup\big(h_{i+1}^{-1}(y)\cap\sigma_{i+1}\big)$, $y\in M$.
According to our construction, we have $h_i\big(h_{i+1}^{-1}(y)\cap
K_i\big)\subset\St(y,\V_i)$. Hence, $h_{i+1}^{-1}(y)\cap K_i$ is
contained in $\Omega(z)$ for some $z\in h_i(K_i)$. Therefore,
$d_0\big(h_{i+1}^{-1}(y)\cap K_{i}\big)<\theta$. Since
$h_{i+1}^{-1}(y)\cap\sigma_{i+1}$ is $0$-dimensional, Lemma 2.1
implies that $d_0(h_{i+1}^{-1}(y)\cap K_{i+1})<\theta$. Obviously the map
$h^*=h_m$ is $\delta/2$-homotopic to $q_2\circ\alpha$ and $h^*\in
C(K,M;\theta)$. This completes the proof of the claim.

To finish the proof of the proposition, we apply Lemma 2.2 with $G$
being the set of all simplicially factorizable maps from $C(X,M)$
and $U(i)=C(X,M;1/i)$, $i\geq 1$ (we can apply Lemma 2.2 because of
the claim). Hence, for every simplicially factorizable map $g\in
C(X,M)$ and a positive number $\e$ there exists a map
$g'\in\bigcap_{i=1}^\infty U(i)$ which is $\e$-homotopic to $g$.
Finally, we observe that $\bigcap_{i=1}^\infty U(i)$ consists of
$0$-dimensional maps.
\end{proof}

\begin{pro}\label{n-1}
Let $M$ be a completely metrizable space with the
$AP(n,0)$-property, $X$ a compactum and $Z$ an $F_{\sigma}$-subset
of $X$ with $\dim Z\leq n-1$. Then there exists a
$G_{\delta}$-subset $\mathcal H\subset C(X,M)$ with the following
properties:
\begin{itemize}
\item $Z$ is contained in the union of trivial components of the
fibers of $g$ for all $g\in\mathcal H$;
\item for every simplicially
factorizable map $g\in C(X,M)$ and every $\e>0$ there exists
$g'\in\mathcal H$ which is $\e$-homotopic to $g$.
\end{itemize}
\end{pro}

\begin{proof}
We represent $Z$ as the union of an increasing sequence
$(Z_i)_{i\geq 1}$ with all $Z_i$ being closed in $Z$. For every
$\e>0$ let
$$\mathcal H(Z_i,\e)=\{g\in C(X,M): F(g,\e)\cap Z_i=\emptyset\},$$ where
$$F(g,\e)=\cup\{C: \mbox{$C$ is a component
of a fiber of $g$ with $\diam C\geq\e$}\}.$$  It is easily seen that
each $\mathcal H(Z_i,\e)$ is open in $C(X,M)$. So, the set $\mathcal
H=\bigcap_{i=1}^\infty\mathcal H(Z_i,1/i)$ is $G_{\delta}$ and $Z$
is contained in the union of trivial components of the fibers of $g$
for every $g\in\mathcal H$.

To prove the second item of our proposition, we first consider the
particular case when $X$ is a polyhedron and $Z$ is a (compact)
subpolyhedron of $X$.

{\em Claim $1$. Suppose that, in addition to hypotheses of
Proposition $2.4$, $X$ is a polyhedron and $Z$ is a
subpolyhedron of $X$. Then, for every $g_0\in C(X,M)$ and $\delta>0$
there exists $g\in\mathcal H=\bigcap_{i=1}^\infty\mathcal H(Z,1/i)$
which is $\delta$-homotopic to $g_0$.}

The proof of Claim 1 is a slight modification of the proof of
\cite[Theorem 2.2]{km}. Let $g_0\in C(X,M)$ and $\delta>0$. We take
an open neighborhood $W$ of $Z$ in $X$, a retraction $r\colon W\to
Z$ and a function $\alpha\colon X\to\I$ such that
$\alpha^{-1}(0)=Z$, $\alpha^{-1}(1)=X\backslash W$. Since $X$ is a
compact $ANR$, we can choose $W$ so small that $g_0|W$ is
$\delta/2$-homotopic to $(g_0\circ r)|W$. Next, denote by $\pi\colon
Z\times\I\to Z$ the projection and consider the map $\varphi\colon
W\to Z\times\I$, $\varphi (x)=(r(x),\alpha(x))$. Obviously,
$(g_0\circ r)|W=(g_0\circ \pi\circ\varphi)|W$. So, $(g_0\circ
\pi\circ\varphi)|W$ is $\delta/2$-homotopic to $g_0|W$. Because $Z$
is a polyhedron, so is $Z\times\I$. Hence, every map from
$C(Z\times\I,M)$ is simplicially factorizable. Moreover, $\dim
Z\times\I\leq n$ and $M\in AP(n,0)$. Therefore, by Proposition 2.3,
there exists a $0$-dimensional map $h\colon Z\times\I\to M$ which is
$\delta/2$-homotopic to the map $g_0\circ\pi$. Then
$(h\circ\varphi)|W$ is $\delta/2$-homotopic to $(g_0\circ
\pi\circ\varphi)|W$. Consequently, $(h\circ\varphi)|W$ is
$\delta$-homotopic to $g_0|W$. By the Homotopy Extension Theorem,
there exits a map $g\in C(X,M)$ such that $g$ is $\delta$-homotopic
to $g_0$ and $g|\overline{U}=(h\circ\varphi)|\overline{U}$, where
$U$ is an open neighborhood of $Z$ in $X$ with $\overline{U}\subset
W$.

To finish the proof of Claim 1, it remains to show that
$g\in\mathcal H$. Let $C$ be a subcontinuum of $g^{-1}(y)$ for some
$y\in M$ and let $Z\cap C\neq\emptyset$. We are going to prove that
$C\subset Z$. Otherwise, there would be a subcontinuum $C'\subset
C\cap U$ such that $C'\cap Z\neq\emptyset$ and $C'\backslash
Z\neq\emptyset$. Then $g(C')=h(\varphi(C'))=y$ and, according to the
definition of $\alpha$, $\varphi (C')$ is a non-degenerate continuum
in $h^{-1}(y)$. Since $h$ is $0$-dimensional, this is a
contradiction. Hence, $C\subset Z$. Using again that $\dim(h)\leq
0$, we conclude that $\varphi(C)$ is a point. On the other hand
$\varphi(C)=C\times\{0\}\subset Z\times\I$. Therefore, $C$ should be
a trivial continuum.

Now, consider the general case of Proposition 2.4.

{\em Claim $2$. Let  $g_0\in C(X,M)$ be a simplicially factorizable
map and $\delta$, $\eta$ positive numbers. Then for any $i$ there
exists a simplicially factorizable map $g\in\mathcal H(Z_i,\eta)$
which is $\delta$-homotopic to $g_0$.}

Since $g_0$ is simplicially factorizable, there exist a finite
simplicial complex $L$ and maps $q_1\colon X\to L$, $q_2\colon L\to
M$ with $g_0=q_2\circ q_1$. Let $\V=\{B_{\varrho}(y,\delta/4): y\in
M\}\in cov(M)$ and $\W=q_2^{-1}(\V)\in cov(L)$. Next, choose a finite cover
$\U\in cov(X)$ of $X$ with $\mesh(\U)<\eta$
 such that:
\begin{itemize}
\item at most $n$ elements of the family $\gamma=\{U\in\U: U\cap
Z_i\neq\emptyset\}$ can have a common point;
\item there exists a map $h\colon\mathcal N(\U)\to L$ such that
$h\circ f_{\U}$ is $\W$-homotopic to $q_1$, where $\mathcal N(\U)$
is the nerve of $\U$ and $f_{\U}\colon X\to\mathcal N(\U)$ is the
natural map.
\end{itemize}
Then $q_2\circ h\circ f_{\U}$ is $\delta/2$-homotopic to $g_0$ and
the subpolyhedron $K$ of $\mathcal N(\U)$ generated by the family
$\gamma$ is of dimension $\leq n-1$. So, according to Claim 1, there
exists a map $g_1\colon\mathcal N(\U)\to M$ such that $g_1$ is
$\delta/2$-homotopic to $q_2\circ h$ and $K$ is contained in the
union of trivial components of the fibers of $g_1$. Since all fibers of $f_{\U}$
are of diameter
$<\eta$ and $f_{\U}(Z_i)\subset K$, the map $g=g_1\circ
f_{\U}$ belongs to $\mathcal H(Z_i,\eta)$. Moreover, $g$ is
obviously simplicially factorizable and $\delta$-homotopic to $g_0$.
This completes the proof of Claim 2.

Finally, the proof of Proposition 2.4 follows from Lemma 2.2 (with
$G$ being the set of all simplicially factorizable maps from
$C(X,M)$ and $U(i)=\mathcal H(Z_i,1/i)$ for every $i\geq 1$) and
Claim 2.
\end{proof}

\section{Parametric Lelek maps}
 In this section we are going to prove Theorem 1.1. Everywhere in
 this section we suppose that $(M,\varrho)$ is a given complete metric space and
$\mu=\{W_\nu:\nu\in\Lambda\}$, $\mu_1=\{G_\nu:\nu\in\Lambda\}$ are
locally finite open covers of $M$ such that

\begin{enumerate}
\item[(*)] $\overline{G}_\nu\subset W_\nu$ and $W_\nu\in AP(n_\nu,0)$ with $0\leq n_\nu$ for
every $\nu\in\Lambda$.
\end{enumerate}

Obviously, Theorem 1.1 follows directly from next theorem.

\begin{thm}
Let $f\colon X\to Y$ be a perfect map between metrizable spaces with
$\dim\triangle(f)\leq m$. Suppose $(M,\varrho)$ is a complete metric
space and $\mu$, $\mu_1$ two locally finite open covers of $M$
satisfying condition $(*)$. Then there is a $G_{\delta}$-set
$\mathcal{H}\subset C(X,M)$ such that any simplicially factorizable
map in $C(X,M)$ can be homotopically approximated by maps from
$\mathcal{H}$ and every $g\in\mathcal{H}$ has the following
property: for any $y\in Y$ and $\nu\in\Lambda$ the restriction
$g|\big(f^{-1}(y)\cap g^{-1}(\overline{G}_\nu)\big)$ is an
$(m-n_{\nu})$-dimensional Lelek map from $f^{-1}(y)\cap
g^{-1}(\overline{G}_\nu)$ into $\overline{G}_\nu$.
\end{thm}

The proof of Theorem 3.1 consists of few propositions. We can assume
that $n_\nu\leq m$ for all $\nu\in\Lambda$. For every $y\in Y$,
$\e>0$, $\nu\in\Lambda$ and $g\in C(X,M)$ we denote by
$F_\nu(g,\e,y)$ the union of all continua $C$ of diameter $\geq\e$
such that $C\subset f^{-1}(y)\cap g^{-1}(z_C)$ for some
$z_C\in\overline{G}_\nu$. Note that each $F_\nu(g,\e,y)$ is compact
as a closed subset of $f^{-1}(y)$. Then, for a fixed $\eta>0$, let
$\mathcal{H}(y,\e,\eta)$ be the set of all $g\in C(X,M)$ such that
$d_{m-n_\nu}\big(F_\nu(g,\e,y)\big)<\eta$ for every $\nu\in\Lambda$.
The last inequality means that $F_\nu(g,\e,y)$ can be covered by an
open family $\gamma$ in $X$ such that $\mesh(\gamma)<\eta$ and no
more $m-n_\nu+1$ elements of $\gamma$ have a common point. If
$V\subset Y$, then $\mathcal{H}(V,\e,\eta)$ denotes the intersection
of all $\mathcal{H}(y,\e,\eta)$, $y\in V$.

\begin{lem}
For every $y\in Y$ and every $g\in \mathcal{H}(y,\e,\eta)$ there
exists a neighborhood $V_y$ of $y$ in $Y$ and $\delta_y>0$ such that
if $y'\in V_y$ and $\varrho(h(x),g(x))<\delta_y$ for all  $x\in
f^{-1}(y')$, then $h\in\mathcal{H}(y',\e,\eta)$.
\end{lem}

\begin{proof}
Since $\mu_1$ is locally finite and $f^{-1}(y)$ is compact, for
every $y\in Y$ and $g\in C(X,M)$ there exists a neighborhood
$O_g(y)$ of $g(f^{-1}(y))$ such that the family
$\Lambda_g(y)=\{\nu\in\Lambda: O_g(y)\cap G_\nu\neq\emptyset\}$ is
finite.

 Assume the lemma is not true for some $y\in Y$ and $g\in
\mathcal{H}(y,\e,\eta)$. Then there exists a local base of
neighborhoods $V_i$ of $y$, points $y_i\in V_i$ and functions
$g_i\in C(X,M)$ such that $g_i|f^{-1}(y_i)$ is $1/i$-close to
$g|f^{-1}(y_i)$ but $g_i\not\in\mathcal{H}(y_i,\e,\eta)$. It is
easily seen that for some $k$ and all $i\geq k$ we have
$g_i(f^{-1}(y_i))\subset O_g(y)$ . Consequently, for every $i\geq k$
there exists $\nu(i)\in\Lambda_g(y)$ such that
$d_{m-n_{\nu(i)}}\big(F_{\nu(i)}(g_i,\e,y_i)\big)\geq\eta$ (because
$g_i\not\in\mathcal{H}(y_i,\e,\eta)$). So, each
$F_{\nu(i)}(g_i,\e,y_i)$ doesn't have any open cover of $\mesh<\eta$
and order $\leq m-n_{\nu(i)}+1$. Since the family $\Lambda_g(y)$ is
finite, there exists $\nu(0)\in\Lambda_g(y)$ with $\nu(i)=\nu(0)$
for infinitely many $i$. With out loss of generality, we may suppose
that $\nu(i)=\nu(0)$ for all $i\geq k$. This implies that
$g_i(f^{-1}(y_i))\cap\overline{G}\neq\emptyset$, $i\geq k$, where
$G=G_{\nu(0)}$. Consequently,
$g(f^{-1}(y))\cap\overline{G}\neq\emptyset$. Since $g\in
\mathcal{H}(y,\e,\eta)$, $F_{\nu(0)}(g,\e,y)$ can be covered by an
open family $\gamma$ in $X$ of order $\leq m-n_{\nu(0)}+1$ and
$\mesh(\gamma)<\eta$. Let $U=\cup\gamma$. To obtain a contradiction,
it suffice to show that $F_{\nu(i)}(g_i,\e,y_i)\subset U$ for some
$i\geq k$. Indeed, otherwise for every $i\geq k$ there would exist points
$x_i\in F_{\nu(i)}(g_i,\e,y_i)\backslash U$,
$z_i\in\overline{G}$, and a component $C_i$ of $f^{-1}(y_i)\cap
g_i^{-1}(z_i)$ containing $x_i$ with $\diam C_i\geq\e$. Using that
$P=f^{-1}(\{y_i\}_{i=k}^\infty\cup\{y\})$ is a compactum, we can
suppose that $\{x_i\}_{i=k}^\infty$ converges to a point $x_0\in
f^{-1}(y)$, $\{z_i\}_{i=k}^\infty$ converges to a point
$z_0\in\overline{G}$ and $\{C_i\}_{i=k}^{\infty}$ (considered as a
sequence in the space of all closed subsets of $P$ equipped with the
Vietoris topology) converges to a closed set $C\subset f^{-1}(y)\cap
g^{-1}(z_0))$. It is easily seen that $C$ is connected and $\diam
C\geq\e$. Hence, $C\subset F_{\nu(0)}(g,\e,y)\subset U$. So, $x_i\in
C_i\subset U$ for some $i$, a contradiction.
\end{proof}

Now, we are in a position to show that the sets
$\mathcal{H}(Y,\e,\eta)$ are open in $C(X,M)$.

\begin{pro}\label{open}
For any closed set $F\subset Y$ and any $\e, \eta>0$, the set
$\mathcal{H}(F,\e,\eta)$ is open in $C(X,M)$.
\end{pro}

\begin{proof}
Let $g_0\in\mathcal{H}(F,\e,\eta)$. Then, by Lemma 3.2, for every
$y\in F$ there exist a neighborhood $V_y$ and a positive
$\delta_y\leq 1$ such that $h\in\mathcal{H}(y',\e,\eta)$ provided
$y'\in V_y$ and $h|f^{-1}(y')$ is $\delta_y$-close to
$g_0|f^{-1}(y')$.  The family $\{V_y\cap Y:y\in F\}$ can be supposed
to be locally finite in $F$. Consider the set-valued lower
semi-continuous map $\varphi\colon F\to (0,1]$,
$\varphi(y)=\cup\{(0,\delta_z]:y\in V_z\}$. By \cite[Theorem 6.2,
p.116]{rs}, $\varphi$ admits a continuous selection $\beta\colon
F\to (0,1]$. Let $\overline{\beta}:Y\to (0,1]$ be a continuous
extension of $\beta$ and $\alpha=\overline{\beta}\circ f$. It
remains only to show that if $g\in C(X,M)$ with
$\varrho\big(g_0(x),g(x)\big)<\alpha(x)$ for all $x\in X$, then
$g\in\mathcal{H}(F,\e,\eta)$. So, we take such a $g$ and fix $y\in
F$. Then there exists $z\in F$ with $y\in V_{z}$ and
$\alpha(x)\leq\delta_{z}$ for all $x\in f^{-1}(y)$. Hence,
$\varrho\big(g(x),g_0(x)\big)<\delta_z$ for each $x\in f^{-1}(y)$.
According to the choice of $V_z$ and $\delta_z$,
$g\in\mathcal{H}(y,\e,\eta)$. Therefore, $\mathcal{H}(F,\e,\eta)$ is
open in $C(X,M)$.
\end{proof}

To prove Theorem 3.1 it suffices to show that if $g\in
C(X,M)$ is a simplicially factorizable map and $\delta\in
C(X,(0,1])$, then for any $\e,\eta>0$ there exists a simplicially
factorizable map $g_{\e\eta}\in\mathcal{H}(Y,\e,\eta)$ which is
$\delta$-homotopic to $g$.  Indeed, since any set
$\mathcal{H}(Y,\e,\eta)$ is open, Lemma 2.2 would imply that every
simplicially factorizable map is homotopically approximated by simplicially
factorizable maps from
$\mathcal{H}=\bigcap_{i,j=1}^{\infty}\mathcal{H}(Y,1/i,1/j)$. But
for every  $g\in\mathcal H$, $y\in Y$, $\nu\in\Lambda$
and $i\ge1$ we have
$d_{m-n_\nu}\big(F_{\nu}(g,1/i,y)\big)=0$.
So, $\dim F_{\nu}(g,1/i,y)\leq m-n_\nu$. Then $\dim
F_\nu(g,y)\leq m-n_\nu$, where
$F_{\nu}(g,y)=\bigcup_{i=1}^{\infty}F_\nu(g,1/i,y)$ . On the other
hand, $F_\nu(g,y)$ is the union of all non-trivial continua
contained in the fibers of $g|g^{-1}(\overline{G}_\nu)\cap
f^{-1}(y)$. Therefore, $\mathcal{H}$ consists of maps $g$ such that
$g|g^{-1}(\overline{G}_\nu)\cap f^{-1}(y)\colon
g^{-1}(\overline{G}_\nu)\cap f^{-1}(y)\to \overline{G}_\nu$ is
$(m-n_\nu)$-dimensional Lelek map for every $\nu\in\Lambda$ and
$y\in Y$.

Next proposition shows that all simplicially factorizable maps in
$C(X,M)$ can be homotopically approximated by simplicially
factorizable maps from $\mathcal{H}(Y,\e,\eta)$ provided the set
$\mathcal{H}_p(L)$ is homotopically dense in $C(N,M)$ , where
$p\colon N\to L$ is any perfect $m$-dimensional $PL$-map between two
simplicial complexes $N, L$ equipped with the $CW$-topology and
$\mathcal{H}_p(L)$ is the set of all maps $q\colon N\to M$ such that
$q|q^{-1}(\overline{G}_\nu)\cap p^{-1}(z)\colon
q^{-1}(\overline{G}_\nu)\cap p^{-1}(z)\to \overline{G}_\nu$ is an
$(m-n_\nu)$-dimensional Lelek map for every $\nu\in\Lambda$ and
$z\in L$. Recall that $p\colon N\to L$ is a $PL$-map (resp., a
simplicial map) if $p$ maps every simplex $\sigma$ of $N$ into
(resp., onto) some simplex of $L$ and $p$ is linear on $\sigma$.

\begin{pro}
Let $X$, $Y$, $f$ and $M$ satisfy the hypotheses of Theorem $3.1$.
Suppose the set $\mathcal{H}_p(L)$ is homotopically dense in
$C(N,M)$ for any perfect $m$-dimensional $PL$-map $p\colon N\to L$.
Then for any simplicially factorizable map $g\in C(X,M)$ and any
$\delta\in C(X,(0,1])$, $\e,\eta>0$ there exists a simplicially
factorizable $h\in\mathcal{H}(Y,\e,\eta)$ such that $h$ is
$\delta$-homotopic to $g$.
\end{pro}

\begin{proof}
For fixed $\delta\in C(X,(0,1])$ and a simplicially factorizable map
$g\in C(X,M)$ we are going to find a simplicially factorizable
$h\in\mathcal{H}(Y,\e,\eta)$ such that
$\varrho(g(x),h(x))<\delta(x)$ for all $x\in X$, where $\e$ and
$\eta$ are arbitrary positive reals. Since $g$ is simplicially
factorizable, there exists a simplicial complex $D$ and maps
$g_D\colon X\to D$, $g^D\colon D\to M$ with $g=g^D\circ g_D$. The
metric $\varrho$ induces a continuous pseudometric $\varrho_D$ on
$D$, $\varrho_D(x,y)=\varrho(g^D(x),g^D(y))$. By \cite{ca} and
\cite{si}, $D$ being a stratifiable $ANR$ is a neighborhood retract
of a locally convex space. Hence, we can apply \cite[Lemma 8.1]{bv}
to find an open cover $\U$ of $X$ satisfying the following
condition: if $\alpha\colon X\to K$ is a $\U$-map into a paracompact
space $K$ (i.e., $\alpha^{-1}(\omega)$ refines $\U$ for some
$\omega\in cov(K)$), then there exists a map $q'\colon G\to D$,
where $G$ is an open neighborhood of $\overline{\alpha(X)}$ in $K$,
such that $g_D$ and $q'\circ\alpha$ are $\delta/2$-homotopic with
respect to the metric $\varrho_D$. Let $\U_1$ be an open cover of
$X$ refining $\U$ with $\mesh\U_1<\min\{\e,\eta\}$ and
$\inf\delta(U)>0$ for all $U\in\U_1$.

Next, according to \cite[Theorem 6]{bv}, there exists an open cover
$\V$ of $Y$ such that: for any $\V$-map $\beta\colon Y\to L$ into a
simplicial complex $L$ we can find an $\U_1$-map $\alpha\colon X\to
K$ into a simplicial complex $K$ and a perfect $m$-dimensional
$PL$-map $p\colon K\to L$ with $\beta\circ f=p\circ\alpha$.
We can assume that $\V$ is locally finite. Take
$L$ to be the nerve of the cover $\V$ and $\beta\colon Y\to L$ the
corresponding natural map. Then there are a simplicial complex $K$
and maps $p$ and $\alpha$ satisfying the above conditions. Hence,
the following diagram is commutative:
$$
\begin{CD}
X@>{\alpha}>>K\cr @V{f}VV @VV{p}V\cr Y@>{\beta}>>L\cr
\end{CD}
$$
Since $K$ is paracompact, the choice of the cover $\U$ guarantees
the existence of a map $q_D\colon G\to D$, where $G\subset K$  is an
open neighborhood of $\overline{\alpha(X)}$, such that $g_D$ and
$h_D=q_D\circ\alpha$ are $\delta/2$-homotopic with respect to
$\varrho_D$. Then, according to the definition of $\varrho_D$,
$h'=g^D\circ q_D\circ\alpha$ is $\delta/2$-homotopic to $g$ with
respect to $\varrho$. Replacing the triangulation of $K$ by a
suitable subdivision, we may additionally assume that no simplex of
$K$ meets both $\overline{\alpha(X)}$ and $K\backslash G$. So, the
union $N$ of all simplexes $\sigma\in K$ with
$\sigma\cap\overline{\alpha(X)}\neq\emptyset$ is a subcomplex of $K$
and $N\subset G$. Moreover, since $N$ is closed in $K$, $p\colon
N\to L$ is a perfect $m$-dimensional $PL$-map. Therefore, we have
the following commutative diagram, where $q=g^D\circ q_D$:

\begin{picture}(120,95)(-100,0)
\put(30,10){$L$} \put(0,30){$Y$} \put(12,28){\vector(3,-2){18}}
\put(14,14){\small $\beta$} \put(1,70){$X$}
\put(5,66){\vector(0,-1){25}} \put(-1,53){\small $f$}
\put(11,73){\vector(1,0){45}} \put(30,77){\small $h'$}
\put(12,68){\vector(3,-2){18}} \put(15,56){\small $\alpha$}
\put(31,50){$N$} \put(35,46){\vector(0,-1){25}}
 \put(37,33){\small $p$}
\put(46,58){\vector(4,3){13}} \put(44,64){\small $q$}
 \put(60,70){$M$}
\end{picture}

Now, we shall construct a continuous function $\delta_1:N\to(0,1]$
with $\delta_1\circ\alpha\leq\delta$. Since $\alpha$ is a
$\U_1$-map, there is an open cover $\V_1$ of $N$ such that the cover
$\alpha^{-1}(\V_1)=\{\alpha^{-1}(V):V\in\V_1\}$ refines $\U_1$.
Because $\inf\delta(U)>0$ for any $U\in\U_1$,
$\inf\delta(\alpha^{-1}(V))>0$ for any $V\in\V_1$. We can assume
that $\V_1$ is locally finite and consider the lower semi-continuous
set-valued map $\varphi\colon N\to (0,1]$ defined by $\varphi
(z)=\cup\{(0,\inf\delta(\alpha^{-1}(V))]:z\in V\in\V_1\}$. Then, by
\cite[Theorem 6.2, p. 116]{rs}, $\varphi$ admits a continuous
selection $\delta_1\colon N\to (0,1]$. Obviously,
$\delta_1(z)\leq\inf\delta(\alpha^{-1}(z))$ for all $z\in N$. Hence,
$\delta_1\circ\alpha\leq\delta$.

Since, according to our assumption, $\mathcal{H}_p(L)$ is
homotopically dense in $C(N,M)$, there exists a map
$q_1\in\mathcal{H}_p(L)$ such that $q_1$ is $\delta_1/2$-homotopic
to  $q$.  Let $h=q_1\circ\alpha$. Then $h$ and $q\circ \alpha$ are
$\delta/2$-homotopic because $\delta_1\circ\alpha\leq\delta$. On the
other hand, $q\circ\alpha=h'$ is $\delta/2$-homotopic to $g$. Hence,
$g$ and $h$ are $\delta$-homotopic. Moreover, $h$ is obviously
simplicially factorizable.

It remains to show that $h\in\mathcal{H}(Y,\e,\eta)$. To this end,
we fix $y\in Y$ and $\nu\in\Lambda$, and consider the set
$F_\nu(h,\e,y)$. Recall that $F_\nu(h,\e,y)$ is the union of all
continua of diameter $\geq\e$ such that $C\subset f^{-1}(y)\cap
h^{-1}(a_C)$ for some $a_C\in\overline{G}_\nu$. For any such
continuum $C$ we have $\alpha (C)\subset p^{-1}(\beta(y))\cap
q_1^{-1}(a_C)$. Since the diameters of all fibers of $\alpha$ are
$<\e$ (recall that $\alpha$ is an $\U_1$-map), $\alpha(C)$ is a
non-trivial continuum in $p^{-1}(\beta(y))\cap q_1^{-1}(a_C)$.
Therefore, $\alpha(C)\subset F_\nu(q_1,\beta(y))$, where
$F_\nu(q_1,\beta(y))$ denotes the union of all non-trivial continua
which are contained in the fibers of the restriction
$q_{y\nu}=q_1|\big(p^{-1}(\beta(y))\cap
q_1^{-1}(\overline{G}_\nu)\big)$. Actually, we proved that
$\alpha\big(F_\nu(h,\e,y)\big)\subset F_\nu(q_1,\beta(y))$. Since
$q_1\in\mathcal{H}_p(L)$, $q_{y\nu}\colon p^{-1}(\beta(y))\cap
q_1^{-1}(\overline{G}_\nu)\to\overline{G}_\nu$ is
$(m-n_\nu)$-dimensional Lelek map. Consequently, $\dim
F_\nu(q_1,\beta(y))\leq m-n_\nu$.  So, there exists an open cover
$\gamma$ of $F_\nu(q_1,\beta(y))$ of order $\leq m-n_\nu+1$ (such a
cover $\gamma$ exists because $F_\nu(q_1,\beta(y))$ is metrizable as
a subset of the metrizable compactum $p^{-1}(\beta(y)$). We can
suppose that $\gamma$ is so small that $\alpha^{-1}(\gamma)$ refines
$\U_1$. But $\mesh(\U_1)<\eta$. Consequently,
$\alpha^{-1}(\gamma)\cap F_\nu(h,\e,y)$ is an open cover of
$F_\nu(h,\e,y)$ of order $\leq m-n_\nu+1$ and mesh $<\eta$. This
means that $d_{m-n\nu}\big(F_\nu(h,\e,y)\big)<\eta$. Therefore, we
found a simplicially factorizable map $h\in\mathcal{H}(Y,\e,\eta)$
which is $\delta$-homotopic to $g$.
\end{proof}

In next two lemmas we suppose that $p\colon N\to L$ is an
$m$-dimensional $PL$-map between finite simplicial complexes. As
everywhere in this section, $(M,\varrho)$ is a complete metric space
possessing two locally finite open covers
$\mu=\{W_\nu:\nu\in\Lambda\}$ and $\mu_1=\{G_\nu:\nu\in\Lambda\}$
such that $\overline{G}_\nu\subset W_\nu$, $W_\nu\in AP(n_\nu,0)$
and $0\leq n_\nu\leq m$ for every $\nu\in\Lambda$. For given
$\e,\eta>0$ and $y\in L$ we denote by $\mathcal{H}_p(y,\e,\eta)$ the
set of $g\in C(N,M)$ such that
$d_{m-n_\nu}\big(F_\nu(g,\e,y)\big)<\eta$ for every $\nu\in\Lambda$.
Here, $F_\nu(g,\e,y)$ is the union of all continua $C\subset
p^{-1}(y)\cap g^{-1}(z_C)$ with $z_C\in\overline{G}_\nu$ and $\diam
C\geq\e$. If $B\subset L$, then $\mathcal{H}_p(B,\e,\eta)$ stands
for the intersection of all $\mathcal{H}_p(y,\e,\eta)$, $y\in B$.

\begin{lem} Let $B\subset L$ be a
subcomplex of $L$ and $g_0\in C(N,M)$ be a map such that
$g_0\in\mathcal{H}_p(B)=\bigcap_{i\geq 1}\mathcal{H}_p(B,1/i,1/i)$.
Then, for every $\e>0$, $q\geq 1$ and $g\in C(N,M)$ with
$g|p^{-1}(B)=g_0|p^{-1}(B)$ there exists a map $g_q\in C(N,M)$
extending $g_0|p^{-1}(B)$ such that $g_q$ is $\e$-homotopic to $g$
and $g_q\in\mathcal{H}_p(y,1/q,1/q)$ for all $y\in L$.
\end{lem}

\begin{proof}
We fix $q$ and $g\in C(N,M)$ with $g|p^{-1}(B)=g_0|p^{-1}(B)$. Then,
by Lemma 3.2, for every $y\in B$ there exists a
neighborhood $V_y$ in $L$ and $\delta_y>0$ such that any $h\in
C(N,M)$ belongs to $\mathcal{H}_p(V_y,1/q,1/q)$ provided
$h|p^{-1}(V_y)$ is $\delta_y$-close to $g|p^{-1}(V_y)$ (we can apply
Lemma 3.2 because $g|p^{-1}(B)=g_0|p^{-1}(B)$ yields
$g\in\mathcal{H}_p(y,1/q,1/q)$ for every $y\in B$). Let
$\{V_{y(i)}\}_{i\leq s}$ be a finite subfamily of $\{V_y:y\in B\}$
covering $B$ and $V=\bigcup_{1\leq i\leq s}V_{y(i)}$.

Since $\mu$ is locally finite in $M$, $g(N)$ meets only finitely
many $W_j=W_{\nu(j)}$, $j=1,..,k$. For any $j\leq k$ let
$P_j=g^{-1}(\overline{G}_{\nu(j)})$ and $U_j^1, U_j^2$ be open
subsets of $N$ such that $P_j\subset
U_j^1\subset\overline{U}_j^1\subset
U_j^2\subset\overline{U}_j^2\subset g^{-1}(W_j)$. We also choose
$\e_0>0$ satisfying the following condition:

\begin{itemize}
\item If $h\in C(N,M)$ with $\varrho\big(g(x),h(x)\big)<\e_0$ for
every $x\in N$, then $\{\nu\in\Lambda:h(N)\cap W_\nu\neq\emptyset\}$
is contained in $\{\nu(j):j=1,..,k\}$, and for all $j$ we have
  $h^{-1}(\overline{G}_{\nu(j)})\subset U_j^1$ and
$h(\overline{U}_j^2)\subset W_j$.
\end{itemize}

Let $\delta_0=\min\{\e,\e_0,\delta_{y(i)}:i\leq s\}$.
Considering suitable subdivisions of $N$ and $L$, we can suppose
that $p$ is a simplicial map and the following conditions hold:

\begin{itemize}
\item Every simplex $\sigma\in L$ intersecting the set $L\backslash V$
does not meet $B$;
\item Every simplex $\tau\in N$ intersecting the set $\overline{U}_j^1$ does not
meet $N\backslash U_j^2$, $j=1,..,k$.
\end{itemize}

Let $L_0=\cup\{\sigma\in L:\sigma\cap L\backslash V\neq\emptyset\}$
and $N_j=\cup\{\tau\in N:\tau\cap\overline{U}_j^1\neq\emptyset\}$,
$1\leq j\leq k$. Obviously $L_0$ is a subpolyhedron of $L$ disjoint
from $B$ and containing $L\backslash V$. Moreover, each $N_j$ is a
subpolyhedron of $N$ such that $\overline{U}_j^1\subset N_j\subset
U_j^2$. Now, for every $j\leq k$ consider the map $p_j=p|N_j\colon
N_j\to p(N_j)$. Since $p$ is $m$-dimensional, by \cite{bp:96} or
\cite{ht:85}, there exists an $(n_{\nu(j)}-1)$-dimensional sigma-compact
set $Z_j=\bigcup_{i\geq 1}Z_{ij}\subset N_j$ such that each $Z_{ij}$
is compact and $\dim p_j^{-1}(y)\backslash Z_j\leq m-n_{\nu(j)}$ for all
$y\in p(N_j)$. Denote
$$T_{ij}=\{h\in C(N,M):F^j(h,1/q)\cap Z_{ij}=\emptyset\},$$ where
$F^j(h,1/q)$ is the union of all continua $C\subset N_j$ of diameter
$\geq 1/q$ which are contained in fibers of the map $h|N_j$. It is
easily seen that $T_{ij}$ are open in $C(N,M)$.

{\em Claim. For any $h\in B_\varrho(g,\delta_0)$, $\eta>0$ and
$i,j\geq 1$ there exists $h_{ij}\in T_{ij}\cap B_\varrho(h,\eta)\cap
B_\varrho(g,\delta_0)$ which is $\eta$-homotopic to $h$.}

Indeed, $h\in B_\varrho(g,\delta_0)$ and $\delta_0\leq\e_0$ imply
that $h(N_j)\subset W_j$ and $h(N)\cap W_\nu=\emptyset$ if
$\nu\neq\nu(j)$ for all $j$. Let
$\eta_1=\min\{\eta,\delta_0-\varrho\big(g(x),h(x)\big):x\in N\}$.
Since $h|N_j$ is simplicially factorizable (as a map whose domain is
a polyhedron), $W_j\in AP(n_{\nu(j)},0)$ and $Z_{ij}$ is a compact
subset of $N_j$ with $\dim Z_{ij}\leq n_{\nu(j)}-1$, according to
Proposition 2.4, there is a map $h'\colon N_j\to W_j$ which is
$\eta_1$-homotopic to $h|N_j$ and the union of all non-trivial
components of the fibers of $h'$ is disjoint from $Z_{ij}$. By the
Homotopy Extension Theorem, $h'$ admits an extension $h_{ij}\in
C(N,M)$ with $h_{ij}$ being $\eta_1$-homotopic to $h$. Obviously,
$h_{ij}\in T_{ij}$ and $h_{ij}\in B_\varrho(g,\delta_0)$.

The above claim allows us to apply Lemma 2.2 for the set
$B_\varrho(g,\delta_0)$ and the sequence $\{T_{ij}\}_{i,j\geq 1}$ to
obtain a map $h_1\in C(N,M)$ such that
$h_1\in\bigcap_{i,j=1}^{\infty}T_{ij}$ and $h_1$ is
$\delta_0$-homotopic to $g$. Let $h_2\colon p^{-1}(L_0\cup
B)\to M$ be defined by $h_2(x)=g_0(x)$ if $x\in p^{-1}(B)$ and
$h_2(x)=h_1(x)$ if $x\in p^{-1}(L_0)$. Obviously, $h_2$ is
$\delta_0$-homotopic to $g$. Since $L_0\cup B$ is a subpolyhedron of
$L$ and $p$ is a simplicial map, $p^{-1}(L_0\cup B)$ is a
subpolyhedron of $N$. Then, by the Homotopy Extension Theorem, there
exists a map $g_q\in C(N,M)$ extending $h_2$ with $g_q$ being
$\delta_0$-homotopic to $g$.

It remains only to show that $g_q\in\mathcal{H}_p(y,1/q,1/q)$ for
all $y\in L$. This is true if $y\in V$. Indeed, then $y$ belongs to
some $V_{y(i)}$. Since $g_q|V_{y(i)}$ is $\delta_0$-close to
$g|V_{y(i)}$ and $\delta_0\leq\delta_{y(i)}$,
$g_q\in\mathcal{H}_p(y,1/q,1/q)$ according to the choice of
$V_{y(i)}$ and $\delta_{y(i)}$.  If $y\in L_0$, then
$g_q|p^{-1}(y)=h_1|p^{-1}(y)$. Since $h_1$ is $\delta_0$-close to
$g$ and $\delta_0\leq\e_0$, $h_1^{-1}(\overline{G}_{\nu(j)})\subset
U_j^1\subset N_j$ for every $j\leq k$. So, $p^{-1}(y)\cap
g_q^{-1}(\overline{G}_{\nu(j)})=p^{-1}(y)\cap
h_1^{-1}(\overline{G}_{\nu(j)})=p_j^{-1}(y)\cap
g_q^{-1}(\overline{G}_{\nu(j)})$, $j\leq k$. On the other hand,
$h_1\in\bigcap_{i,j=1}^{\infty}T_{ij}$ implies that every
restriction $h_1|h_1^{-1}(\overline{G}_{\nu(j)})\colon
h_1^{-1}(\overline{G}_{\nu(j)})\to\overline{G}_{\nu(j)}$ has the
following property: the union of all non-trivial components of the
fibers of $h_1|h_1^{-1}(\overline{G}_{\nu(j)})$ is contained in
$N_j\backslash Z_j$. Hence, the union of all non-trivial components
of the fibers of $g_q|p^{-1}(y)\cap g_q^{-1}(\overline{G}_{\nu(j)})$
is contained in $p_j^{-1}(y)\backslash Z_j$. Since $\dim
p_j^{-1}(y)\backslash Z_j\leq m-n_{\nu(j)}$, every
$g_q|p^{-1}(y)\cap g_q^{-1}(\overline{G}_{\nu(j)})$ is an
$(m-n_{\nu(j)})$-dimensional Lelek map. But $g_q(N)$ doesn't meet
any $\overline{G}_\nu$ except for $\nu\in\{\nu(j):j=1,..,k\}$.
Therefore, $g_q\in\mathcal{H}_p(y,1/i,1/i)$ for every $i\geq 1$.
\end{proof}

\begin{lem}
Let $B\subset L$ be a subcomplex of $L$ and $g_0\in C(N,M)$ be a map
such that $g_0\in\mathcal{H}_p(B)=\bigcap_{i\geq
1}\mathcal{H}_p(B,1/i,1/i)$. Then, for every $\delta>0$ there exists
$\overline{g}_0\in\mathcal{H}_p(L)$ which is $\delta$-homotopic to
$g_0$ and $\overline{g}_0|p^{-1}(B)=g_0|p^{-1}(B)$
\end{lem}

\begin{proof}
Each set $\mathcal{H}_p(L,1/i,1/i)$ is open in $C(N,M)$ according to
Proposition 3.3. So, by Lemma 3.5, we can apply Lemma 2.2 (with
$U(i)$ being in our case $\mathcal{H}_p(L,1/i,1/i)$ and
$A=p^{-1}(B)$) to find a map $\overline{g}_0\in\bigcap_{i\geq
1}\mathcal{H}_p(L,1/i,1/i)$ which is $\delta$-homotopic to $g_0$ and
$\overline{g}_0|p^{-1}(B)=g_0|p^{-1}(B)$. Finally,
$\overline{g}_0\in\mathcal{H}_p(L)$ because
$\mathcal{H}_p(L)=\bigcap_{i\geq 1}\mathcal{H}_p(L,1/i,1/i)$.
\end{proof}

Next proposition completes the proof of Theorem 3.1. We suppose that
$p\colon N\to L$ is a perfect $m$-dimensional $PL$-map between
simplicial complexes, $(M,\varrho)$ a complete metric space and
$\mu=\{W_\nu:\nu\in\Lambda\}$, $\mu_1=\{G_\nu:\nu\in\Lambda\}$ are
locally finite open covers of $M$ with $\overline{G}_\nu\subset
W_\nu$ and $W_\nu\in AP(n_\nu,0)$ for every $\nu\in\Lambda$, where
$n_\nu\leq m$ are integers. If $B\subset A\subset L$, denote by
$\mathcal{H}_p(A,B)$ the set of the maps $g\in C(p^{-1}(A),M)$ such
that $g|g^{-1}(\overline{G}_\nu)\cap p^{-1}(y)\colon
g^{-1}(\overline{G}_\nu)\cap p^{-1}(y)\to \overline{G}_\nu$ is
$(m-n_\nu)$-dimensional Lelek map for every $\nu\in\Lambda$ and
$y\in B$. When $A=B$, we write $\mathcal{H}_p(B)$ instead of
$\mathcal{H}_p(A,B)$.

\begin{pro}
Let $p\colon N\to L$ and $(M,\varrho)$ be as above. Then the set
$\mathcal{H}_p(L)$ is homotopically dense in $C(N,M)$.
\end{pro}

\begin{proof}
As usual, the simplicial complexes $N$ and $L$ are equipped with the
$CW$-topology. But when consider a diameter of any subset of $N$ we
mean the diameter with respect to the standard metric generating the
metric topology of $N$. According to the notations in this section,
for every sets $B\subset A\subset L$ and $\e,\eta>0$, let
$\mathcal{H}_p(A,B,\e,\eta)$ be the set of all $g\in C(p^{-1}(A),M)$
such that any $d_{m-n_\nu}\big(F_\nu(g,\e,y)\big)<\eta$ for all
$y\in B$ and $\nu\in\Lambda$. Although the domain of $g$ is the set
$A$ (not the whole space $N$), we use the same notation
$F_\nu(g,\e,y)$ to denote the union of all continua $C\subset
p^{-1}(y)\cap g^{-1}(z_C)$ with $z_C\in\overline{G}_\nu$ and $\diam
C\geq\e$. Let us also denote by $\mathcal{H}_p(A,B)$ the set of the
maps $g\in C(p^{-1}(A),M)$ such that $g|g^{-1}(\overline{G}_\nu)\cap
p^{-1}(y)\colon g^{-1}(\overline{G}_\nu)\cap p^{-1}(y)\to
\overline{G}_\nu$ is an $(m-n_\nu)$-dimensional Lelek map for every
$\nu\in\Lambda$ and $y\in B$. This means that $\dim F_\nu(g,y)\leq
m-n_\nu$ for any $y\in B$ and $\nu\in\Lambda$, where
$F_\nu(g,y)=\bigcup_{i=1}^{\infty}F_\nu(g,1/i,y)$. It is easily seen
that $\mathcal{H}_p(A,B)=\bigcap_{i\geq
1}\mathcal{H}_{p}(A,B,1/i,1/i)$.

Now, let us finish the proof of Proposition 3.7. Fix $g\in C(N,M)$
and $\delta\in C(N,(0,1])$. We are going to find
$h\in\mathcal{H}_p(L)$ which is $\delta$-homotopic to $g$. To this
end, let $L^{(i)}$, $i\geq 0$, denote the $i$-dimensional skeleton
of $L$ and $L^{(-1)}=\emptyset$. We put $h_{-1}=g$ and construct
inductively a sequence $(h_i:N\to M)_{i\geq 0}$ of maps such that
\begin{itemize}
\item $h_{i}|p^{-1}(L^{(i-1)})=h_{i-1}|p^{-1}(L^{(i-1)})$;
\item $\displaystyle h_{i}$ is $\displaystyle\frac{\delta}{2^{i+2}}$-homotopic to $h_{i-1}$;
\item $h_i\in\mathcal{H}_p(L,L^{(i)})$.
\end{itemize}

Assuming that the map $h_{i-1}:N\to M$ has been constructed,
consider the complement $L^{(i)}\setminus L^{(i-1)}=\sqcup_{j\in
J_i}\overset{\circ}\sigma_j$, which is the discrete union of open
$i$-dimensional simplexes.  Since $h_{i-1}|\sigma_j$ belongs to
$\mathcal{H}_p(\sigma_j,\sigma^{(i-1)}_j)$ for any simplex
$\sigma_j\in L^{(i)}$, we can apply Lemma 3.6 to find a map
$g_j:p^{-1}(\sigma_j)\to M$ such that
\begin{itemize}
\item $g_j$ coincides with $h_{i-1}$ on the set $p^{-1}(\sigma^{(i-1)}_j)$;
\item $g_j$ is $\displaystyle\frac{\delta}{2^{i+2}}$-homotopic to
$h_{i-1}$;
\item $g_j\in\mathcal{H}_p(\sigma_j,\sigma_j)$.
\end{itemize}

Next, define a map $q_i:p^{-1}(L^{(i)})\to M$ by the formula
$$q_i(x)=\begin{cases} h_{i-1}(x)&\mbox{if $x\in
p^{-1}(L^{(i-1)})$;}\\ g_j(x)&\mbox{if $x\in p^{-1}(\sigma_j)$.}
\end{cases}$$ It can be shown that $q_i$ is
$\displaystyle\frac{\delta}{2^{i+2}}$-homotopic to
$h_{i-1}|p^{-1}(L^{(i)})$. Since $p^{-1}(L^{(i)})$ is a
subpolyhedron of $N$, we can apply the Homotopy Extension Theorem to
find a continuous extension $h_i:N\to M$ of the map $q_i$ which is
$\displaystyle\frac{\delta}{2^{i+2}}$-homotopic to $h_{i-1}$.
Moreover, $h_i\in\mathcal{H}_p(L,L^{(i)})$ because
$h_{i-1}\in\mathcal{H}_p(L,L^{(i-1)})$ and
$g_j\in\mathcal{H}_p(\sigma_j,\sigma_j)$ for any $j$. This completes
the inductive step.

Then the limit map $h=\lim_{i\to\infty}h_i:N\to M$ is well-defined,
continuous and $\delta$-homotopic to $g$ (the last two properties of
$h$ hold because $h$ has this properties for any simplex from $N$
and because of the definition of the $CW$-topology on $N$). Finally,
since $h|p^{-1}(L^{(i)})=h_i|p^{-1}(L^{(i)})$ and $h_i\in\mathcal{H}_p(L,L^{(i)})$
for every $i$,
$h\in\mathcal{H}_p(L)$.
\end{proof}

\section{Almost $AE(n,0)$-spaces}

We already observed that if $M$ is an $\LC[n]$-space, then $M\in
AP(n,0)$ if and only if $M$ has the following property:

\begin{itemize}
\item for every map $g\in C(\I^n,M)$ and every $\e>0$ there exists a
$0$-dimensional map $g'\in C(\I^n,M)$ which is $\e$-close to $g$.
\end{itemize}
Any space having the above property will be referred as
{\em almost $AP(n,0)$}. Obviously, every $\LC[n-1]$ almost $AP(n,0)$-space has the
$AP(n-1,0)$-property.

We are going to establish an analogue of Theorem 3.1 for almost
$AP(n,0)$-spaces.

\begin{lem}\label{anr}
Every complete $\LC[n-1]$-space $M$ admits a complete metric $\varrho$
generating its topology and satisfying the following condition: If
$Z$ is an $n$-dimensional space, $A\subset Z$ its closed set and
$h\colon Z\to M$, then
for every function $\alpha:Z\to(0,1]$ and every map $g\colon A\to M$
with $\varrho(g(z),h(z))<\alpha(z)/8$ for all $z\in A$ there exists
a map $\bar{g}\colon Z\to M$ extending $g$ such that
$\varrho(\bar{g}(z),h(z))<\alpha(z)$ for all $z\in Z$.
\end{lem}

\begin{proof}
We embed $M$ in a Banach space $E$ as a closed subset. Since the
Hilbert cube is the image of the $n$-dimensional Menger compactum
under an $n$-invertible map \cite{ad:84}, we can find a metric
space $E(n)$ with $\dim E(n)\leq n$ and a perfect $n$-invertible
surjection $p\colon E(n)\to E$. Here, $p$ is $n$-invertible means
that every map from at most $n$-dimensional space into $E$
can be lifted to a map into $E(n)$.

Since $M\in\LC[n-1]$, there exist a neighborhood $W$ of $M$ in $E$
and a map $q\colon p^{-1}(W)\to M$ extending the restriction
$p|p^{-1}(M)$. For every open $U\subset M$ let $T(U)=W\backslash
p(q^{-1}(M\backslash U))$. Obviously, $T(U)\subset W$ is open,
$T(U)\cap M=U$ and $q(p^{-1}(T(U)))=U$. Let $\mathcal{T}$ be the
collection of all pairs $(U,V)$ of open sets in $M$ such that
$\overline{conv(V)}\subset T(U)$, where $\overline{conv(V)}$ is the
closed convex hull of $V$ in $E$. Now, consider the family
$$\mathcal{T}=\{(U,V):U,V\hbox{~}\mbox{are open
in}\hbox{~}M\hbox{~}\mbox{and}\hbox{~}\overline{conv(V)}\subset
T(U)\}.$$ The family $\mathcal{T}$ has the following properties: (i)
for any $z\in M$ and its neighborhood $U$ in $M$ there is a
neighborhood $V\subset U$ of $z$ with $(U,V)\in\mathcal{T}$; (ii)
for any $(U,V)\in\mathcal{T}$ and open sets $U',V'\subset M$ we have
$(U',V')\in\mathcal{T}$ provided $U\subset U'$ and $V'\subset V$. By
\cite[Proposition 2.3]{dm:97}, there exists a complete metric
$\varrho$ on $M$ such that for every $z\in M$ and $r\in (0,1)$ the
pair of open balls $(B_{\varrho}(z,r),(B_{\varrho}(z,r/8))$ belongs
to $\mathcal{T}$.

Suppose we are given an $n$-dimensional space $Z$, its closed subset
$A\subset Z$ and two maps $h\colon Z\to M$ and
$g\colon A\to M$ such that $\varrho(g(z),h(z))<\alpha(z)/8$ for all $z\in A$, where
$\alpha\in C(Z,(0,1])$.
Consider the
set-valued map $\phi\colon Z\to E$, $\phi(z)=g(z)$ if $z\in A$ and
$\phi(z)=\overline{conv(B_{\varrho}(h(z),\alpha(z)/8))}$ if
$z\not\in A$. Then $\phi$ is lower semi-continuous and has closed
and convex values in $E$. So, by the Michael convex-valued selection
theorem \cite{em:56}, $\phi$ has a continuous selection $g_1$. Next,
we lift $g_1$ to a map $g_2\colon Z\to E(n)$.  According to the
definition of $\mathcal{T}$, every
$p^{-1}\big(\overline{conv(B_{\varrho}(h(z),\alpha(z)/8))}\big)$ is
contained in $q^{-1}\big(B_{\varrho}(h(z),\alpha(z))\big)$. Hence,
$\overline{g}=q\circ g_2$ is the required extension of $g$.
\end{proof}

We also need the following lemma whose proof is similar to that
one of Lemma 2.2.

\begin{lem}
Let $G\subset C(X,M)$, where $(M,\varrho)$ is a complete metric space.
Suppose $(U(i)_{i\geq 1}$ is a sequence of open subsets of $C(X,M)$
such that
\begin{itemize}
\item for any $g\in G$, $i\geq 1$ and any function $\eta\in
C(X,(0,1])$ there exists $g_i\in B_\varrho(g,\eta)\cap U(i)\cap G$
which is $\eta$-close to $g$.
\end{itemize}
Then, for any $g\in G$ and $\e\colon X\to (0,1]$ there exists
$g'\in\bigcap_{i=1}^{\infty} U(i)$ which is and $\e$-close to $g$.
If, in addition, all $g_i$ can be chosen such that $g_i|A=g_0|A$ for
some $g_0\in C(X,M)$ and $A\subset X$, then $g'|A=g_0|A$.
\end{lem}

Next two propositions are analogues of Proposition 2.3 and Proposition 2.4.

\begin{pro}\label{n-1}
Let $M$ be a complete $\LC[n-1]$ almost $AP(n,0)$-space. Then for every
$n$-dimensional compactum $X$ there exists a dense
$G_{\delta}$-subset $\mathcal H\subset C(X,M)$ of $0$-dimensional maps.
\end{pro}

\begin{proof}
First of all, let us note that since $\dim X\leq n$ and $M$ is $\LC[n-1]$, the set of all
simplicially factorizable maps from $C(X,M)$ is dense in $C(X,M)$.
Analyzing  the proof of Proposition 2.3 and using Lemma 4.2 instead of Lemma 2.2,
one can see that it suffices to establish the following claim:

{\em Claim. If $K$ is a finite $n$-dimensional polyhedron,
then every map $g\colon K\to M$ can be approximated by a
$0$-dimensional map $g'\colon K\to M$.}

Since the sets $C(K,M;\eta)$ consisting of maps $g\in C(K,M)$ with
$d_0(g^{-1}(g(x))<\eta$, $\eta>0$, are open in $C(K,M)$, according
to Lemma 4.2, it is enough to show that for every $g\in C(K,M)$ and
$\eta>0$ there exists $g'\in C(K,M;\eta)$ which is $\eta$-close to
$g$. To this end, we equipped $M$ with a complete metric satisfying
the hypotheses of Lemma 4.1 and fix $0<\eta\leq 1$ and $g\in
C(K,M)$. Since $M\in\LC[n-1]$ and $M$ is almost $AP(n,0)$, $M\in
AP(n-1,0)$. So, there exists a $0$-dimensional map $h\colon
K^{(n-1)}\to M$ which is $\eta/16$-close to $g|K^{(n-1)}$. Here,
$K^{(n-1)}$ is the $(n-1)$-dimensional skeleton of $K$. By Lemma
4.1, $h$ can be extended to a map $\bar{h}\in C(K,M)$ with
$\varrho(\bar{h}(x),g(x))<\eta/2$ for all $x\in K$. Let $\sigma_j$,
$j=1,..,k$, be all $n$-dimensional simplexes of $K$. Then
$K\backslash K^{(n-1)}$ is a disjoint union of the open simplexes
$\overset{\circ}\sigma_j$. Let
$K_i=K^{(n-1)}\bigcup_{j=1}^{j=i}\sigma_j$, $i=1,..,k$ We are going
to construct by induction maps $h_i\colon K\to M$, $i=0,..,k$,
satisfying the following conditions:
\begin{itemize}
\item $h_0=\bar{h}$;
\item $h_i|K_i$ belongs to $C(K_i,M;\eta)$, $1\leq i\leq k$;
\item $h_i$ and $h_{i+1}$ are $\displaystyle (\eta/2k)$-close, $i=0,..,k-1$.
\end{itemize}

Assume
that $h_i$ has  already been constructed. Since $h_i|K_i$ belongs to
$C(K_i,M;\eta)$, every fiber $h_i^{-1}(y)\cap K_i$ of $h_i|K_i$,
$y\in h_i(K_1)$, is covered by a finite open and disjoint family
$\Omega(y)$ in $K_i$ with $\mesh(\Omega(y))<\eta$. Using that
$h_i(K_i)$ is compact, we find $\displaystyle 0<\delta_i<\eta/2k$
such that if $\varrho(y,h_i(K_i))<\delta_i$ for some $y\in M$, then there is $z\in
h_i(K_i)$ with $h_i^{-1}(B_\varrho(y,\delta_i))\cap K_i\subset\Omega(z)$.
Since $M$ is almost $AP(n,0)$, there exists a $0$-dimensional
map $p_i\colon\sigma_{i+1}\to M$ which is $\delta_i/8$-close to
$h_i|\sigma_{i+1}$. By Lemma 4.1, $p_i$ can be
extended to a map $h_{i+1}\colon K\to M$ being $\delta_i$-close to
$h_i$. To show
that $h_{i+1}|K_{i+1}$ belongs to $C(K_{i+1},M;\eta)$, we observe
that $h_{i+1}^{-1}(y)\cap K_{i+1}=\big(h_{i+1}^{-1}(y)\cap
K_i\big)\cup\big(h_{i+1}^{-1}(y)\cap\sigma_{i+1}\big)$, $y\in M$.
According to our construction, we have $h_i\big(h_{i+1}^{-1}(y)\cap
K_i\big)\subset B_\varrho(y,\delta_i)\cap h_i(K_i)$. Hence,
$h_{i+1}^{-1}(y)\cap K_i\subset h_i^{-1}(B_\varrho(y,\delta_i))\cap K_i
\subset\Omega(z)$ for some $z\in h_i(K_i)$. Therefore,
$d_0\big(h_{i+1}^{-1}(y)\cap K_{i}\big)<\eta$. Since
$h_{i+1}^{-1}(y)\cap\sigma_{i+1}$ is $0$-dimensional, Lemma 2.1
implies that $d_0(h_{i+1}^{-1}(y)\cap K_{i+1})<\eta$. Obviously $h_k\in
C(K,M;\eta)$ and
$h_k$ is $\eta/2$-close to $\bar{h}$. Hence, $g'=h_k$ is $\eta$-close to $g$.
This completes the proof of the claim.
\end{proof}

\begin{pro}
Let $M$ be a complete $\LC[n-1]$ almost $AP(n,0)$-space. Then for every
$n$-dimensional compactum $X$ and its $F_{\sigma}$-subset $Z$ with $\dim Z\leq n-1$
there exists a dense
$G_{\delta}$-subset $\mathcal H\subset C(X,M)$ of maps $g$
such that $Z$ is contained in the union of trivial components of the fibers of $g$.
\end{pro}

\begin{proof}
Following the proof of Proposition 2.4 and using Lemma 4.2 and Proposition 4.3 instead of
Lemma 2.2 and Proposition 2.3, respectively, it suffices to prove the following analogue of Claim 1
from Proposition 2.4.

{\em Claim. Suppose $X$ is a finite $n$-dimensional polyhedron and $Z$ a subpolyhedron of $X$ with
$\dim Z\leq n-1$. Then for every $g_0\in C(X,M)$ and $\delta>0$ there exists
$g\in\mathcal{H}=\bigcap_{i=1}^{\infty}\mathcal{H}(Z,1/i)$ which is $\delta$-close to $g_0$.}

Here, $\mathcal{H}(Z,1/i)$ is the set of all $g\in C(X,M)$ with $F(g,1/i)\cap Z=\emptyset$.
We follow the proof of Claim 1, Proposition 2.4 and use the same notations. The first difference is that
we take $W$ to be a neighborhood of $Z$ in $X$ such that $g_0|W$ and $(g_0\circ r)|W$ are
$\delta/16$-close. Then, $(g_0\circ\pi\circ\varphi)|W$ is $\delta/16$-close to $g_0|W$. Next,
we use by Proposition 4.3 to choose a $0$-dimensional map $h\colon Z\times\I\to M$
which is $\delta/16$-close to $g_0\circ\pi$. So, $(h\circ\varphi)|W$ is $\delta/8$-close to $g_0|W$.
Finally, take a neighborhood $U$ of $Z$ in $X$ with $\overline{U}\subset W$, and use
Lemma 4.1 to find an extension $g\in C(X,M)$ of
$(h\circ\varphi)|\overline{U}$ with $g$ being  $\delta$-close to $g_0$. Then $g\in\mathcal{H}$.
\end{proof}

We establish now an analogue of Theorem 3.1 for almost
$AP(n,0)$-spaces.

\begin{thm}
Let $f\colon X\to Y$ be a perfect map between metrizable spaces with
$\dim X\leq n$ and $M$ is a complete $\LC[n-1]$-space. Suppose
$\mu=\{W_\nu:\nu\in\Lambda\}$ and $\mu_1=\{G_\nu:\nu\in\Lambda\}$
are locally finite open covers of $M$ such that
$\overline{G}_\nu\subset W_\nu$ and each $W_\nu$ is almost
$AP(m_\nu,0)$ for every $\nu\in\Lambda$. Then there is a dense
$G_{\delta}$-set $\mathcal{H}\subset C(X,M)$  of maps $g$ such that
for any $y\in Y$ and $\nu\in\Lambda$ the restriction
$g|\big(f^{-1}(y)\cap g^{-1}(\overline{G}_\nu)\big)$ is an
$(n-m_{\nu})$-dimensional Lelek map.
\end{thm}

\begin{proof}
Following the notations and the proof of Theorem 3.1, we can assume
that $m_\nu\leq n$ for all $\nu$. Let $\mathcal{H}(y,\e,\eta)$
denote the set of all $g\in C(X,M)$ such that
$d_{n-m_\nu}\big(F_\nu(g,\e,y)\big)<\eta$ for every $\nu\in\Lambda$,
where $y\in Y$ and $\e,\eta>0$ are fixed. Moreover, for any
$F\subset Y$, let $\mathcal{H}(F,\e,\eta)$ be the intersection of
all $\mathcal{H}(y,\e,\eta)$, $y\in F$. One can establish a lemma
analogical to Lemma 3.2. Then, as in Proposition 3.3, we can show
that any $\mathcal{H}(F,\e,\eta)$ is open in $C(X,M)$, where
$F\subset Y$ is closed.

Observe that $\dim\triangle(f)\leq m$ because $X$, as a space of
dimension $\leq n$, admits a uniformly $0$-dimensional map into
$\I^n$, see \cite{mk}. Moreover, we can assume that the simplicial
complex $K$ in Proposition 3.4 is $n$-dimensional. So, we can apply
Proposition 3.4 in the present situation to show that any
$\mathcal{H}(Y,\e,\eta)$ is dense in $C(X,M)$ provided for any
perfect $PL$ map $p\colon N\to L$ with $\dim N\leq n$ the set
$\mathcal{H}_p(L)$ is dense in $C(N,M)$. Here, $\mathcal{H}_p(L)$ is
the set of all maps $q\colon N\to M$ such that
$q|q^{-1}(\overline{G}_\nu)\cap p^{-1}(z)$ is an
$(n-m_\nu)$-dimensional Lelek map for every $\nu\in\Lambda$ and
$z\in L$. Hence, it remains to show that $\mathcal{H}_p(L)$ is dense
in $C(N,M)$ for any $PL$-map $p\colon N\to L$ with $\dim N\leq n$.
And this follows from the proof of Lemma 3.5, Lemma 3.6 and
Proposition 3.7 with the only difference that now we replace the
application of Lemma 2.2, Proposition 2.3 and Proposition 2.4 by
Lemma 4.2, Proposition 4.3 and Proposition 4.4, respectively.
\end{proof}

\begin{cor}
Let $M$ be a complete $\LC[n-1]$-space such that each point $z\in M$
has a neighborhood which is almost $AP(n,0)$. Then, for every
perfect map $f\colon X\to Y$ between metric spaces with $\dim X\leq
n$, there is a dense $G_{\delta}$-set $\mathcal{H}\subset C(X,M)$
consisting of maps $g$ such that every restriction $g|f^{-1}(y)$,
$y\in Y$, is a $0$-dimensional map.
\end{cor}

Any manifold $M$ modeled on the $n$-dimensional Menger cube has the
$AP(n,0)$-property, see Corollary 6.5 below. So, Theorem 1.1 holds
for such a space $M$. But Theorem 1.1 does not provide any
information about the density of the set $\mathcal{H}$ in $C(X,M)$
except that every simplicially factorizable map in $C(X,M)$ can be
approximated by maps from $\mathcal{H}$. Next proposition shows
that, in this special case, the set $\mathcal{H}$ is dense in
$C(X,M)$ with respect to the uniform convergence topology.

\begin{pro}
Let $f\colon X\to Y$ be a perfect map between metrizable spaces with
$\dim\triangle(f)\leq m$ and $M$ be a manifold modeled on the
$n$-dimensional Menger cube. Then there is a $G_{\delta}$-set
$\mathcal{H}\subset C(X,M)$ dense in $C(X,M)$ with respect to the
uniform convergence topology such that for any $g\in\mathcal{H}$ and
$y\in Y$ the restriction $g|f^{-1}(y)$ is an $(m-n)$-dimensional
Lelek map.
\end{pro}

\begin{proof}
Let $\mathcal{H}$ be the set from the proof of Theorem 3.1. To show
that $\mathcal{H}$ is dense in $C(X,M)$ equipped with the uniform
convergence topology, we used an idea from the proof of
\cite[Corollary 2.8]{km}. According to \cite{be}, for any $\e>0$
there exists an $n$-dimensional polyhedron $P\subset M$ of piecewise
embedding dimension $n$ and maps $u\colon M\to P$ and $v\colon P\to
M$ such that $u$ is a retraction, $v$ is $0$-dimensional, $v\circ u$
is $\e/2$-close to the identity $id_M$. Since every $ANR$ of
piecewise embedding dimension $n$ has the $AP(n,0)$-property (see
\cite[Propsition 2.1]{km}), according to Theorem 3.1, for every
$g\in C(X,M)$ there is $g'\colon X\to P$ such that $g'$ is
$\delta$-close to $u\circ g$ and $g'|f^{-1}(y)$ is an
$(m-n)$-dimensional Lelek map for all $y\in Y$. Here $\delta>0$ is
chosen such that $dist(v(x),v(y))<\e/2$ for any $x,y\in P$ which are
$\delta$-close. Then $v\circ g'$ is $\e$-close to $g$ and since $v$
is $0$-dimensional, $v\circ g'\in\mathcal{H}$.
\end{proof}

\section{Properties of $AP(n,0)$-spaces}
In this section we investigate the class of $AP(n,0)$-spaces.

\begin{lem}
Let $K$ be a polyhedron of dimension $\leq n$ and $L\subset K$ a
subpolyhedron. Suppose $(X,\varrho)$ is a complete metric space
possessing the $AP(n,0)$-property and $g_0\in C(K,X)$ with $\dim
g_0^{-1}(g_0(x))\leq 0$ for all $x\in L$. Then for every $\delta>0$
there exists a $0$-dimensional map $g\colon K\to X$ which is
$\delta$-homotopic to $g_0$ and $g|L=g_0|L$.
\end{lem}

\begin{proof}
We already observed in Section 2 that all sets $C(K,X,\e)$,
$\e>0$, consisting of maps $h\in C(K,X)$ with $d_0(h^{-1}(h(x)))<\e$
for every $x\in K$ are open in $C(K,X)$. Since every map from
$C(K,X)$ is simplicially factorizable, $C(K,X,\e)$ are homotopically dense in
$C(K,X)$ according to Proposition 2.3.

{\em Claim $1$. For every $x\in L$ and $j\geq 1$ there exist
$\e_x>0$ and a neighborhood $U_x$ of $x$ in $K$ satisfying the
following condition: If  $h\in C(K,X)$ and $Z\subset K$ with
$\varrho\big(g_0(y),h(y)\big)<\e_x$ for all $y\in Z\cup\overline{U}_x$, then
$d_0(h^{-1}(h(y))\cap Z)< 1/j$ for any $y\in\overline{U}_x$.}

The proof of this claim is similar to the proof of Lemma 3.2.

\smallskip
{\em Claim $2$. Let $h\in C(K,X)$ with $h|L=g_0|L$. Then, for every
$\eta>0$ and $j\geq 1$ there exists $h_j\in C(K,X,1/j)$ such that
$h_j|L=g_0|L$ and $h_j$ is $\eta$-homotopic to $h$.}

We fix $j$ and $\eta>0$. Choose finitely many points $x(i)\in L$,
$i\leq k$, positive reals $\e_{x(i)}$ and neighborhoods $U(x_i)$ in
$K$ satisfying the hypotheses of Claim 1 such that $L\subset
U=\bigcup_{i\leq k}U(x_i)$. Taking a smaller neighborhood, if
necessarily, we can suppose that
$\varrho\big(g_0(x),h(x)\big)<\eta_1/2$ for all $x\in\overline{U}$,
where $\eta_1=\min\{\eta,\e_{x(i)}:i\leq k\}$. Consider a
triangulation $T$ of $K$ such that $\sigma\in T$ and $\sigma\cap
K\backslash U\neq\emptyset$ imply $\sigma\cap L=\emptyset$. Now, let
$N$ be the subpolyhedron of $K$ given by $N=\cup\{\sigma\in
T:\sigma\cap K\backslash U\neq\emptyset\}$. Obviously, $N$ and $L$
are disjoint. Since $X\in AP(n,0)$ and $\dim N\leq n$, by
Proposition 2.3, there exists a $0$-dimensional map $g_N\in C(N,X)$
which is $\eta_1/2$-homotopic to $h|N$ (we can apply Proposition 2.3
because $h|N$ is simplicially factorizable as a map with a
polyhedral domain). The map $h'\colon N\cup L\to X$, $h'|N=g_N$ and
$h'|L=g_0|L$, is $\eta_1/2$-homotopic to $h|(N\cup L)$. So, by the
Homotopy Extension Theorem, $h'$ admits an extension $h_j\colon K\to
X$ which is $\eta_1/2$-homotopic to $h$. It remains only to show
that $d_0(h_j^{-1}(h_j(x)))< 1/j$ for any $x\in K$. To this end,
observe that $\varrho\big(g_0(x),h_j(x)\big)<\eta_1$ for all $x\in
\overline{U}$. Because $K=N\cup\overline{U}$, for every $x\in K$ we
have $h_j^{-1}(h_j(x))=\big(h_j^{-1}(h_j(x))\cap
N\big)\cup\big(h_j^{-1}(h_j(x))\cap\overline{U}\big)$. Since
$h_j|N=g_N$ and $g_N$ is $0$-dimensional, $\dim h_j^{-1}(h_j(x))\cap
N\leq 0$. On the other hand,
$h_j^{-1}(h_j(x))\cap\overline{U}=h_j^{-1}(h_j(y))\cap\overline{U}$
for some $y\in\overline{U}$. So, there exists $m\leq k$ with
$y\in\overline{U}_{x(m)}$. Since $h_j|\overline{U}$ is
$\eta_1$-close to $g_0|\overline{U}$,
$\varrho\big(g_0(z),h_j(z)\big)<\e_{x(m)}$ for all $z\in
\overline{U}$. Hence, according to the choice of $U_{x(m)}$ and
$\e_{x(m)}$, $d_0\big(h_j^{-1}(h_j(y))\cap\overline{U})\big)<1/j$.
Therefore, $d_0\big(h_j^{-1}(h_j(x))\cap\overline{U})\big)<1/j$.
Finally, by Lemma 2.1, $d_0\big(h_j^{-1}(h_j(x))\big)<1/j$. This
completes the proof of Claim 2.

We are in a position to complete the proof of Lemma 5.1. Because of
Claim 2, we can apply Lemma 2.2 (with $G$ being in our case the set
$\{h\in C(K,X):h|L=g_0|L\}$ and $U(j)=C(K,X,1/j)$) to obtain a map
$g\in\bigcap_{j=1}^{\infty}C(K,X,1/j)$ such that $g|L=g_0|L$ and $g$
is $\delta$-homotopic to $g_0$. Obviously,
$g\in\bigcap_{j=1}^{\infty}C(K,X,1/j)$ yields $\dim g\leq 0$.
\end{proof}

\begin{pro} For any space $X$ we have:
\begin{enumerate}
\item If $X$ has the $AP(n,0)$-property, then every open
subset of $X$ also has this property.
\item If $X$ is completely metrizable, then  $X$ has the $AP(n,0)$-property if and only if
it admits a cover by open subsets with that property.
\end{enumerate}
\end{pro}

\begin{proof}
To prove the first item, suppose $U$ is an open subset of $X\in
AP(n,0)$, $\U\in cov(U)$ and $g\in C(\I^n,U)$. Let
$\U'=\U\cup\{X\backslash g(\I^n)\}$. Since $X\in AP(n,0)$, there is
a $0$-dimensional map $g'\in C(\I^n,X)$ such that $g'$ is
$\U'$-homotopic to $g$. Then $g'(\I^n)\subset U$ and there exists a
$\U$-homotopy $h\colon\I^n\to U$ joining $g$ and $g'$.

According to Michael's theorem on local properties \cite{mi}, the
second item will be established if we show that: (i) A space has the
$AP(n,0)$-property provided $X$ is a discrete sum of spaces with the
same property; (ii) A completely metrizable space has the
$AP(n,0)$-property provided it is a union of two open subspaces with
this property. Condition (i) trivially follows from the definition.
Let us check condition (ii).

Suppose $X$ is a completely metrizable space and $X=X_0\cup X_1$ is
the union of two open subspaces $X_1,X_2\subset X$ with the
$AP(n,0)$-property. Fix an open cover $\U$ of $X$ and a map
$g:\mathbb I^n\to X$ and choose two open sets $W_1,W_2\subset X$
such that $X=W_1\cup W_2$ and $W_i\subset \overline{W}_i\subset X_i$
for $i\in\{1,2\}$. Next, find a complete metric $\varrho$ on $X$
such that
\begin{itemize}
\item $\varrho(X\setminus W_1,X\setminus W_2)\geq 1$;
\item $B_\varrho(\overline{W}_i,1)\subset X_i$ for $i\in\{1,2\}$ and
\item each set of diameter $<1$ in $(X,\varrho)$ lies in some $U\in\U$.
\end{itemize}

Let $V=V_1\cap V_2$, where $V_i=B_\rho(\overline{W}_i,1/2)$,
$i=1,2$, and choose a triangulation $T$ of $\mathbb I^n$ such that
for any simplex $\sigma\in T$ we have $\diam g(\sigma)<1/4$. Now,
consider the polyhedra
$$\begin{aligned} K_i&=\cup\{\sigma\in T:g(\sigma)\cap
\overline{W}_i\neq\emptyset\},\hbox{} i\in\{1,2\}
\end{aligned}$$ and

$$\begin{aligned} L_2&=\cup\{\sigma\in T:g(\sigma)\cap
\overline{V}_2\neq\emptyset\}.
\end{aligned}$$

Obviously, $g(K_i)\subset V_i\subset\overline{V}_i\subset X_i$,
$i=1,2$. So, $g(K_0)\subset V$, where $K_0=K_1\cap K_2$. Moreover,
we have
$$K_2\subset g^{-1}(V_2)\subset\overline{g^{-1}(V_2)}\subset
L_2\subset g^{-1}(X_2)$$ and
$$K_0\subset g^{-1}(V)\subset L_2.$$
Choose now positive $\delta\leq\min\{\varrho\big(g(L_2),X\backslash
X_2\big),\varrho\big(g(K_1),X\backslash X_1\big),1/2\}$ such that
$$h^{-1}(h(K_0))\subset g^{-1}(V)$$ for any $h\in C(\I^n,X)$ which
is $\delta$-close to $g$.

Since $X_2\in AP(n,0)$ and $L_2$ is a polyhedron of dimension $\leq
n$, by Proposition 2.3, there exists a $0$-dimensional map
$g_2\colon L_2\to X_2$ which is $\delta$-homotopic to $g|L_2$. Next,
by the Homotopy Extension Theorem, $g_2$ can be extended to a map
$\overline{g}_2\in C(\I^n,X)$ $\delta$-homotopic to $g$. According
to the choice of $\delta$, $\overline{g}_2(L_2)\subset X_2$,
$\overline{g}_2(K_1)\subset X_1$ and
$(\overline{g}_2)^{-1}(\overline{g}_2(K_0))\subset g^{-1}(V)\subset
L_2$. Hence, the restriction $g_0=\overline{g}_2|K_1$ is a
map from $K_1$ into $X_1$ and, for every $x\in K_0$, we have
$g_0^{-1}(g_0(x))\subset g_2^{-1}(g_2(x))$.
The last inclusion implies that
$\dim g_0^{-1}(g_0(x))\leq 0$ for all $x\in K_0$ because $g_2$ is
$0$-dimensional. Since $X_1\in AP(n,0)$, we can apply Lemma 5.1 (for
the polyhedra $K_0\subset K_1$ and the map $g_0$) to obtain a
$0$-dimensional map $g_1\in C(K_1,X_1)$ such that $g_1|K_0=g_0$ and
$g_1$ is $1/2$-homotopic to $g_0$. Finally, consider the map
$g_{12}\in C(\I^n,X)$ defined by $g_{12}|K_1=g_1$ and
$g_{12}|K_2=g_2|K_2$. Since both $g_1$ and $g_2$ are
$0$-dimensional, so is $g_{12}$. Moreover, $g_{12}$ is $1$-homotopic
to $g$ because $\delta\leq 1/2$. So, according to the choice of
$\varrho$, $g_{12}$ is $\U$-homotopic to $g$.
\end{proof}

\begin{pro} If $X$ and $Y$ are complete metric spaces such
that $X\in AP(n,0)$ and $Y\in AP(m,0)$, then $X\times Y\in
AP(n+m,0)$.
\end{pro}

\begin{proof}
Let $\varrho_X$ and $\varrho_Y$ be complete metrics on $X$ and $Y$,
respectively. We fix $g\in C(\I^{n+m},X\times Y)$ and consider the
complete metric $\varrho=\max\{\varrho_X,\varrho_Y\}$ on $X\times
Y$. It suffices to show that for every $\varepsilon>0$ there exists
a $0$-dimensional map $h\in C(\I^{n+m},X\times Y)$ which is
$\varepsilon$-homotopic to $g$. Let $g_X=\pi_X\circ g$ and
$g_Y=\pi_Y\circ g$, where $\pi_X$ and $\pi_Y$ are the projections
from $X\times Y$ onto $X$ and $Y$, respectively. We represent
$\I^{n+m}$ as the union $\bigcup_{i=1}^{n+m+1}A_i$ with each $A_i$
being a $0$-dimensional $G_{\delta}$-subset of $\I^{n+m}$. Then
$Z_X=\I^{n+m}\backslash\bigcup_{i=n+1}^{n+m+1}A_i$ is an
$F_{\sigma}$-subset of $\I^{n+m}$ which is contained in
$\bigcup_{i=1}^{n}A_i$. So, $\dim Z_X\leq n-1$. Since $X\in AE(n,0)$
and every map from $C(\I^{n+m},X)$ is simplicially factorizable,
Proposition 2.4 yields the existence of a  map $h_X\in
C(\I^{n+m},X)$ such that $h_X$ is $\varepsilon$-homotopic to $g_X$
and $F(h_X)=\bigcup_{j=1}^\infty F(h_X,1/j)\subset\I^{n+m}\backslash
Z_X$. Obviously, $F(h_X)$ is an $F_{\sigma}$-set in $\I^{n+m}$ with
$ F(h_X)\subset\bigcup_{i=n+1}^{n+m+1}A_i$. Hence,
$Z_Y=F(h_X)\backslash A_{n+m+1}$ is also an $F_{\sigma}$-set in
$\I^{n+m}$ with $\dim F_Y\leq m-1$. Now, since $Y\in AE(m,0)$, we
may apply again Proposition 2.4 to obtain a map $h_Y\in
C(\I^{n+m},Y)$ which is $\varepsilon$-homotopic to $g_Y$ and
$F(h_Y)=\bigcup_{j=1}^\infty F(h_Y,1/j)\subset\I^{n+m}\backslash
Z_Y$. Then the diagonal product $h=h_X\triangle h_Y\colon\I^{n+m}\to
X\times Y$ is $\varepsilon$-homotopic to $g$. It remains only to
show that $h$ is $0$-dimensional. If $C$ is a non-trivial component
of a fiber of $h$, then $C\subset F(h_X)\cap F(h_Y)\subset
F(h_X)\backslash Z_Y\subset A_{m+n+1}$. Since $A_{m+n+1}$ is
0-dimensional, $C$ should be a point. Therefore, all components of
the fibers of $h$ are trivial, i.e. $\dim h=0$.
\end{proof}

Finally, we are going to show that every arc-wise connected
$AP(n,0)$-compactum is a continuum $(V^n)$ in the sense of P.
Alexandroff. Recall that a compact metric space $(M,\varrho)$ is a
$(V^n)$-continuum \cite{ps} if for any pair of disjoint closed
subsets $A$ and $B$ of $M$ both having non-empty interiors there
exists $\e>0$ such that $d_{n-2}(C)>\e$ for every partition $C$ in
$M$ between $A$ and $B$. It is easily seen that this is a
topological property, i.e., it doesn't depend on the metric
$\varrho$. Obviously, every continuum $(V^n)$ is a Cantor
$n$-manifold (a compactum which is not disconnected by any
$(n-2)$-dimensional closed subset). Moreover, any $(V^n)$ continuum
has a stronger property \cite{hv}: it cannot be decomposed into a
countable union of proper closed subsets $F_i$ with $\dim F_i\cap
F_j\leq n-2$. The compacta with the last property are called strong
Cantor $n$-manifolds \cite{nh}.

\begin{pro}
Every path-connected compactum $M\in AP(n,0)$ is a continuum $(V^n)$.
\end{pro}

\begin{proof}
Let $A$ and $B$ be disjoint closed subsets of $M$ with non-empty
interior and $\varrho$ be a metric on $M$. Since $M$ is
path-connected, we can choose a map $g\colon\I^n\to M$ such that
$g(\I^n)\cap Int(A)\neq\emptyset\neq g(\I^n)\cap Int(B)$. Then there
exists a $0$-dimensional map $g_1\colon\I^n\to M$ which is so close
to $g$ that $g_1(\I^n)$ meets both $Int(A)$ and $Int(B)$. Thus,
$A_1=g_1^{-1}(A)$ and $B_1=g_1^{-1}(B)$ are disjoint closed subsets
of $\I^n$ with non-empty interiors. Since $\I^n$ is a continuum
$(V^n)$ \cite{ps}, there exists $\e>0$ such that $d_{n-2}(C_1)>\e$
for every partition $C_1\subset\I^n$ between $A_1$ and $B_1$.
Because $g_1$ is $0$-dimensional, every $y\in M$ has a neighborhood
$W_y$ such that $g_1^{-1}(W_y)$ splits into a finite disjoint family
of open subsets of $\I^n$ each of diameter $<\e$. Let $\delta$ be
the Lebesgue number of the cover $\{W_y:y\in M\}$. Then
$d_{n-2}(C)\geq\delta$ for any partition $C\subset M$ between $A$
and $B$. Indeed, otherwise there would be a partition $C$ and an
open family $\gamma_C$ in $M$ of order $\leq n+1$ such that
$\gamma_C$ covers $C$ and $\diam (W)<\delta$ for every
$W\in\gamma_C$. Hence, $g_1^{-1}(\gamma_C)$ is an open family in
$\I^n$ of order $\leq n+1$ and covering $g_1^{-1}(C)$. Moreover,
each $g_1^{-1}(W)$, $W\in\gamma_C$, splits into a finite disjoint
open family consisting of sets with diameter $<\e$. Therefore,
$d_{n-2}(g_1^{-1}(C))<\e$. This is a contradiction because
$g_1^{-1}(C)$ is a partition in $\I^n$ between $A_1$ and $B_1$.
\end{proof}

\section{$AP(n,0)$-spaces and general position properties}

The parametric general position properties were introduces in
\cite{bv}. We say that a space $M$ has the
$m\mbox{-}\overline{DD}^{\{n,k\}}$-property, where $m,n,k\geq 0$ are
integers or $\infty$, if for any open cover $\U$ of $M$ and any two
maps $f:\mathbb I^m\times \mathbb I^n\to M$, $g:\mathbb I^m\times
\mathbb I^k\to M$ there exist maps $f':\mathbb I^m\times \mathbb
I^n\to M$, $g':\mathbb I^m\times \mathbb I^k\to M$ which are
$\U$-homotopic to $f$ and $g$, respectively, and $f'(\{z\}\times
\mathbb I^n)\times g'(\{z\}\times \mathbb I^k\}=\emptyset$ for all
$z\in \mathbb I^m$. It is clear that this is exactly the well known
disjoint $n$-disks property provided $m=0$, $n=k$ and $M$ is
$\LC[n]$. When $m=0$, we simply write $\overline{DD}^{\{n,k\}}$
instead of $0\mbox{-}\overline{DD}^{\{n,k\}}$.

\begin{lem}
Let $M$ be completely metrizable having the
$\overline{DD}^{\{n,k\}}$-property. Suppose $X$ is a compactum and
$A,B\subset X$ closed disjoint subsets  with $\dim A\leq n$ and
$\dim B\leq k$. Then every simplicially factorizable map $g\colon
X\to M$ can be homotopically approximated by maps $g'\in C(X,M)$
such that $g'(A)\cap g'(B)=\emptyset$.
\end{lem}

\begin{proof}
Let $g\in C(X,M)$ be simplicially factorizable and $\delta>0$. As in
Claim 2 from the proof of Proposition 2.4, we can find a finite open
cover $\U$ of $X$ such that:
\begin{itemize}
\item at most $n+1$ elements of the family $\gamma_A=\{U\in\U: U\cap
A\neq\emptyset\}$ intersect;
\item at most $k+1$ elements of the family $\gamma_B=\{U\in\U: U\cap
B\neq\emptyset\}$ intersect;
\item $\cup\gamma_A\cap \cup\gamma_B=\emptyset$;
\item there exists a map $h\colon\mathcal{N}(\U)\to M$ such that
$h\circ f_{\U}$ is $\delta/2$-homotopic to $g$, where
$\mathcal{N}(\U)$ is the nerve of $\U$ and $f_{\U}\colon
X\to\mathcal{N}(\U)$ is the natural map.
\end{itemize}

Let $K_A$ and $K_B$ be the subpolyhedra of $\mathcal{N}(\U)$
generated by the families $\gamma_A$ and $\gamma_B$, respectively.
Then $\dim K_A\leq n$, $\dim K_B\leq k$ and $K_A\cap K_B=\emptyset$.
For any simplexes $\sigma\in K_A$ and $\tau\in K_B$ let
$$G(\sigma,\tau)=\{p\in C(\mathcal{N}(\U),M): p(\sigma)\cap
p(\tau)=\emptyset\}.$$ Obviously, each $G(\sigma,\tau)$ is open
in $C(\mathcal{N}(\U),M)$. Using $M\in\overline{DD}^{\{n,k\}}$ and the
Homotopy Extension Theorem, one can show that all $G(\sigma,\tau)$
are homotopically dense in $C(\mathcal{N}(\U),M)$. So, by Lemma 2.2,
there exists $h'\in\bigcap\{G(\sigma,\tau):\sigma\in K_A,\tau\in
K_B\}$ which is $\delta/2$-homotopic to $h$. Then $g'=h'\circ
f_{\U}$ is $\delta$-homotopic to $g$ and $h'(A)\cap
h'(B)=\emptyset$.
\end{proof}

\begin{pro}For a space $X$ we have:
\begin{enumerate}
\item If $X$ is $\LC[0]$, then $X$ has the
$\overline{DD}^{\{0,0\}}$-property if and only if $X$ has no isolated
point;
\item If $X$ is completely metrizable and $X\in\overline{DD}^{\{0,0\}}$,
then $X$ has the $AP(1,0)$-property.
\end{enumerate}
\end{pro}

\begin{proof}
For the first item, see \cite[Proposiion 7]{bv}. To proof the second
item, we fix countably many closed $0$-dimensional subsets $P_i$ of $\I$
consisting of irrational numbers such that
$\dim\I\backslash\bigcup_{i\geq 1}P_i=0$. For every rational $t\in
Q$ and $i\geq 1$ let $U_i(t)=\{g\in C(\I,M): g(t)\not\in g(P_i)\}$.
By Lemma 6.1, each $U_i(t)$ is homotopically dense in $C(\I,M)$. On the other
hand, obviously all $U_i(t)$ are open in $C(\I,M)$. Therefore, by Lemma 2.2, the
set $U=\bigcap\{U_i(t):t\in Q, i\geq 1\}$ is homotopically dense in $C(\I,X)$.
Moreover, for every $g\in U$ and $x\in\I$ we have the following: if
$g^{-1}(g(x))\cap Q\neq\emptyset$, then
$g^{-1}(g(x))\subset\I\backslash\bigcup_{i\geq 1}P_i$; if
$g^{-1}(g(x))\cap Q=\emptyset$, then
$g^{-1}(g(x))\subset\I\backslash Q$. Hence, $U$ consists of
$0$-dimensional maps.
\end{proof}
Next corollary of Propositions 5.3 and 6.2 provides more examples of
$AP(n,0)$-spaces.

\begin{cor}
If $M_i$, $i=1,..n$, are completely metrizable $\LC[0]$-spaces
without isolated points, then $\prod_{i=1}^{i=n}M_i$ has the
$AP(n,0)$-property.
\end{cor}

Here is a characterization of the $\overline{DD}^{\{n,k\}}$-property.

\begin{pro}
A completely metrizable space $M$ has the
$\overline{DD}^{\{n,k\}}$-property, where $n\leq k\leq\infty$, if and
only if $M$ satisfies the following condition $(n,k)$:
\begin{itemize}
\item If $X$ is a compactum and
$A\subset B\subset X$ its $\sigma$-compact subsets with $\dim A\leq n$
and $\dim B\leq k$, then any
simplicially factorizable map $g\in C(X,M)$ can be homotopically
approximated by maps $h\in C(X,M)$ such that $h^{-1}(h(x))\cap B=x$ for
all $x\in A$.
\end{itemize}
\end{pro}

\begin{proof}
Suppose $M\in\overline{DD}^{\{n,k\}}$ and let $A\subset B\subset X$
be two $\sigma$-compact subsets of a compactum $X$ with $\dim A\leq
n$ and $\dim B\leq k$. Then $A=\bigcup_{p\geq 1}A_p$ and
$B=\bigcup_{m\geq 1}B_m$, where $A_p$ and $B_m$ are compact sets of
dimension $\dim A_p\leq n$ and $\dim B_m\leq k$. For every $p,i\geq
1$ let $\omega_i(p)=\{A_{pj}(i):j=1,2,..,s(p,i)\}$ be a family of
closed subsets of $A_p$ such that $\omega_i(p)$ covers $A_p$ and
$\mesh(\omega_i(p))<1/i$. We also choose sequences
$\{B_{mq}(p,i,j)\}_{q=1}^{\infty}$ of closed sets
$B_{mq}(p,i,j)\subset B_m$ with $B_m\backslash
A_{pj}(i)=\bigcup_{q=1}^{\infty}B_{mq}(p,i,j)$, where $p,i\geq 1$
and $j=1,2,..,s(p,i)$. Then the sets
$$G(p,i,j,m,q)=\{g\in C(X,M): g(A_{pj}(i))\cap g(B_{mq}(p,i,j))=\emptyset\}$$
are open in $C(X,M)$ and the intersection $G$ of all $G(p,i,j,m,q)$
consists of maps $g$ with $g^{-1}(g(x))\cap B=x$ for any $x\in A$. Hence,
according to Lemma 2.2, every simplicially factorizable map is
homotopically approximated by maps from $G$ provided the following
is true: For any simplicially factorizable map $g\in C(X,M)$, $\e>0$
and $(p,i,j,m,q)$ with $p,i,m,q\geq 1$ and $1\leq j\leq s(i,p)$, there
exists a simplicially factorizable map $g(p,i,j,m,q)\in G(p,i,j,m,q)$
which is $\e$-homotopic to $g$. Since $M\in\overline{DD}^{\{n,k\}}$ and any
couple $A_{pj}(i)$, $B_{mq}(p,i,j)$ consists of disjoint closed sets in $X$
of dimension $\leq n$ and $\leq k$, respectively,
the last statement follows from Lemma 6.1. Hence $M$ satisfies
condition $(n,k)$.

Suppose now that $M$ satisfies condition $(n,k)$. To show that
$M\in\overline{DD}^{\{n,k\}}$, let $f\colon\I^n\to M$ and
$g\colon\I^k\to M$ be two maps. If $k<\infty$, we denote by $X$ the
disjoint sum $\I^n\biguplus\I^k$ and consider the map $h\colon X\to
M$, $h|\I^n=f$ and $h|\I^k=g$. Since $h$ is simplicially
factorizable (as a map with a polyhedral domain) and $M$ has the
$(n,k)$-property, for every $\e>0$ there is a map $h_1\in C(X,M)$
such that $h_1^{-1}(h_1(x))=x$ for all $x\in\I^n$ and $h_1$ is
$\e$-homotopic to $h$. Then the maps $f_1=h_1|\I^n$ and
$g_1=h_1|\I^k$ are $\e$-homotopic to $f$ and $g$, respectively,
and $f_1(\I^n)\cap g_1(\I^k)=\emptyset$. So,
$M\in\overline{DD}^{\{n,k\}}$. If $n<k=\infty$, we choose $k(1)<\infty$
such that $n\leq k(1)$ and the map $g'=g\circ r_{k(1)}$ is
$\e/2$-homotopic to $g$, where
$r_{k(1)}\colon\I^{\infty}\to\I^{k(1)}$ is the retraction of
$\I^{\infty}$ onto $\I^{k(1)}\subset\I^{\infty}$ defined by
$r_{k(1)}((x_1,x_2,..))=(x_1,x_2,..,x_{k(1)},0,0,..)$. Then the map
$h\colon\I^n\biguplus\I^{\infty}\to M$, $h|\I^n=f$ and
$h|\I^{\infty}=g'$, is simplicially factorizable. Hence, as in the
previous case, we can use the $(n,\infty)$-property of $M$ to show
that $M\in\overline{DD}^{\{n,\infty\}}$. If $n=k=\infty$, we
homotopically approximate both $f$ and $g$ by maps $f'=f\circ
r_{n(1)}$ and $g'=g\circ r_{k(1)}$, respectively, and proceed as in
the first case.
\end{proof}

\begin{cor}
Every completely metrizable space $\displaystyle
M\in\overline{DD}^{\{n-1,n-1\}}$ has the $AP(n,0)$-property. In
particular, every manifold modeled on the $n$-dimensional Menger
cube or the $n$-dimensional N\"{o}beling space has the
$AP(n,0)$-property.
\end{cor}

\begin{proof}
Let $\I^n=A\cup B$, where $B=\I^n\backslash A$ is $0$-dimensional
and $A=\bigcup_{i\geq 1}A_i$ is $\sigma$-compact with $\{A_i\}$
being a sequence of closed $(n-1)$-dimensional subsets of $\I^n$.
Since $C(\I^n,M)$ consists of simplicially factorizable maps, by
Proposition 6.4, every map in $C(\I^n,M)$ is homotopically
approximated by maps $g\in C(\I^n,M)$ such that $g^{-1}(g(x))\cap
A=x$ for all $x\in A$. Let us show that any such $g$ is
$0$-dimensional. Indeed, since for every $x\in\I^n$ the intersection
$g^{-1}(g(x))\cap A$ can have at most one point,
$g^{-1}(g(x))=(g^{-1}(g(x))\cap A)\cup (g^{-1}(g(x))\cap B)$ is
$0$-dimensional.

The second part of the corollary follows from the fact that any
$n$-dimensional Menger manifold, as well as any manifold modeled on
the  $n$-dimensional N\"{o}beling space, is a complete
$\LC[n-1]$-space with the disjoint $n$-disks property \cite{be}. So,
any such a space has the $\overline{DD}^{\{n-1,n-1\}}$-property.
\end{proof}

The following notions were introduced in \cite{bv} and \cite{bck}.
Let A be a closed subset of a space $X$.  We say that $A$ is:
\begin{itemize}
\item a (homotopical) $Z_n$-set in $X$ if for any
an open cover $\U$ of $X$  and a map $f:\II^n\to X$ there is a map
$g:\II^n\to X$ such that $g(\II^n)\cap A=\emptyset$ and $g$ is
 $\U$-near ($\U$-homotopic) to $f$;
\item homological $Z_n$-set in $X$
if $H_k(U,U\setminus A)=0$ for all open sets $U\subset
X$ and all $k<n+1$.
\end{itemize}
Here, $H_{*}(U,U\setminus A)=0$ are the relative singular homology
groups with integer coefficients. Each homotopical $Z_n$-set in $X$
is both a $Z_n$-set and a homological $Z_n$-set in $X$. The converse
is not always true.

\smallskip
\begin{thm}\label{t1} Let $X$ be a locally path-connected $AP(n,0)$-space.
Then every $x\in X$ is a homological $Z_{n-1}$-point in $X$.
\end{thm}

\begin{proof} By \cite[Corollary 8.4]{bck}, it suffices to check that
$(x,0)$ is a homotopical $Z_n$-point in $X\times[-1,1]$ for every
$x\in X$. Given a map $f:\I^n\to X\times[-1,1]$ we are going to
homotopically approximate $f$ by a map $g:\I^n\to X\times[-1,1]$
with $(x,0)\notin g(\I^n)$. Let $f=(f_1,f_2)$ where $f_1:\I^n\to X$
and $f_2:\I^n\to[-1,1]$ are the components of $f$. Since $X\in
AP(n,0)$, there is a 0-dimensional map $g_1:\I^n\to X$ that
approximates $f_1$. Because $Z=g^{-1}_1(x)$ is a $0$-dimensional
subset of $\I^n$ and $\{0\}$ is a $Z_0$-set in $[-1,1]$, the map
$f_2$ can be approximated by a map $g_2:\I^n\to[-1,1]$ such that
$0\notin f_2(Z)$. Then for the map $g=(g_1,g_2):\I^n\to
X\times[-1,1]$ we have $(x,0)\notin g(\I^n)$.
\end{proof}

\begin{pro}\label{p1} If $X\in AP(n,0)$ and $Y\in \infty\mbox{-}\overline{DD}^{\{0,0\}}$,
then $X\times Y\in \overline{DD}^{\{n,n\}}$.
\end{pro}

\begin{proof} Let $f=(f_X,f_Y):\I^n\oplus\I^n\to X\times Y$ be a given map.
Since $X\in AP(n,0)$, we  homotopically
approximate the first component $f_X$ by a 0-dimensional map
$g_X:\I^n\oplus\I^n\to X$. Next, use
$\infty\mbox{-}\overline{DD}^{\{0,0\}}$-property of $Y$ to
homotopically approximate the second component $f_Y$ by a map
$g_Y:\I^n\oplus\I^n\to Y$ that is injective on the fibers of $g_X$.
Then the map $g=(g_X,g_Y):\I^n\oplus\I^n\to X\times Y$ is injective,
witnessing that $X\times Y$ has the
$\overline{DD}^{\{n,n\}}$-property.
\end{proof}

We say that a space $X$ has the $AP(\infty,0)$-property if $X\in
AP(n,0)$ for every $n\geq 1$.

\begin{cor}\label{c1} If $X,Y$ are two locally path-connected spaces both
possessing the $AP(\infty,0)$-property and $Y$ contains a dense set
of homotopical $Z_2$-points, then $X\times Y$ has the
$\overline{DD}^{\{\infty,\infty\}}$-property.
\end{cor}

\begin{proof} By Theorem~\ref{t1}, each point of $Y$ is a homological $Z_\infty$-point.
Taking into account that $Y$ has a dense set of homotopical
$Z_2$-points and applying Theorem 26(4,5) of \cite{bv}, we conclude
that the space $Y\in\infty\mbox{-}\overline{DD}^{\{0,0\}}$. Finally,
by Proposition~\ref{p1}, $X\times Y\in
\overline{DD}^{\{\infty,\infty\}}$.
\end{proof}

\begin{cor}\label{c2} If $X,Y$ are two compact $AR$'s both possessing
the $AP(\infty,0)$-property, then $X\times Y\times\I$ is
homeomorphic to the Hilbert cube.
\end{cor}

\begin{proof} Theorem~\ref{t1} implies that every $y\in Y$ is a
homological $Z_\infty$-point in $Y$. Then, according to
\cite[Theorem 26(5)]{bv}, $Y\in\overline{DD}^{\{1,1\}}$.
Consequently, by \cite[Theorem 10(5)]{bv}, each point of $Y$ is a
homotopical $Z_1$-point in $Y$. Now, applying the Multiplication
Formula for homological $Z$-sets \cite[Theorem 17(1)]{bv}), we
conclude that each point of $Y\times\I$ is a homotopical
$Z_2$-point. Therefore, Corollary~\ref{c1} yields that the product
$X\times (Y\times\I)$ has the
$\overline{DD}^{\{\infty,\infty\}}$-property. Being a compact $AR$,
this product is homeomorphic to the Hilbert cube according to the
Toru\'nczyk characterization theorem of $Q$-manifolds.
\end{proof}

In light of the preceding corollary, it is interesting to remark
that there exists a compact space $X\in AR\cap AP(\infty,0)$ which
is not homeomorphic to the Hilbert cube. To present such an example,
we first prove that the $AP(n,0)$-property  is preserved by a
special type of maps.

A map $\pi:X\to Y$ is called an elementary cell-like map if:
\begin{enumerate}
\item $\pi$ is a fine homotopy equivalence, i.e., for every open cover
$\U$ of $Y$ there exists a map $s:Y\to X$ with $\pi\circ s$ being
$\U$-homotopic to the identity of $Y$;
\item the non-degeneracy set $N_\pi=\{y\in Y:|\pi^{-1}(y)|\ne 1\}$ is at most countable;
\item each fiber $\pi^{-1}(y)$, $y\in N_\pi$, is an arc.
\end{enumerate}

\begin{pro} Let $\pi:X\to Y$ be an elementary cell-like map between
metric spaces such that $X$ is complete
and every $x\in \pi^{-1}(N_\pi)$ is a homotopical $Z_n$-point in
$X$. Then $Y\in AP(n,0)$ provided $X\in AP(n,0)$.
\end{pro}

\begin{proof}
Suppose $X\in AP(n,0)$. To prove that $Y\in AP(n,0)$, fix an open
cover $\gamma$ of $Y$ and a map $f:\I^n\to Y$. Let $\gamma_1$ be an
open cover of $Y$ which is star-refinement of $\gamma$. Since $\pi$
is a fine homotopy equivalence, there is a map $s:Y\to X$ such that
$\pi\circ s$ is $\gamma_1$-homotopic to the identity map of $Y$. Let
$D=\{x_i:i\geq 1\}\subset f^{-1}(N_\pi)$ be a sequence such that
$f^{-1}(y)\backslash D$ is $0$-dimensional for every $y\in N_f$.
Then the sets $W_i=\{g\in C(\I^n,X): x_i\notin g(\I^n)\}$ are open
in $C(\I^n,X)$ equipped with the uniform convergence topology.
Moreover,  each $W_i$ is homotopically dense in $C(\I^n,X)$ because
$x_i$ is a homotopical $Z_n$-point in $X$. We also consider the sets
$U_i=\{g\in C(\I^n,X): d_0(g^{-1}(x))<1/i\hbox{~}\mbox{for
all}\hbox{~}x\in X\}$. Since $X\in AP(n,0)$ and all maps from
$C(\I^n,X)$ are simplicially factorizable, it follows from
Proposition 2.3 that any $U_i$ is open and homotopically dense in
$C(\I^n,X)$. This easily implies that the intersections $V_i=W_i\cap
U_i$ are open and homotopically dense in $C(\I^n,X)$. Consequently,
by Lemma 2.2 (with $G$ being $C(\I^n,X)$), there exists a map
$g\in\bigcap_{i=1}^{\infty}V_i$ which is
$\pi^{-1}(\gamma_1)$-homotopic to $s\circ f$. Obviously, $g$ is
$0$-dimensional, $g(\I^n)\subset X\backslash D$ and $f_1=\pi\circ
g:I^n\to Y$ is $\gamma_1$-homotopic to $\pi\circ s\circ f$. Then
$f_1$ is $\gamma$-homotopic to $f$ because $\pi\circ s\circ f$ is
$\gamma_1$-homotopic to $f$. Hence, $Y\in AP(n,0)$.
\end{proof}

Singh \cite{sin} constructed an elementary cell-like map $f:Q\to X$
from the Hilbert cube $Q$ onto a compact $X\in AR$ such that
$X\times\I$ is homeomorphic to $Q$ but $X$ contains no proper
$ANR$-subspace of dimension $\ge 2$. By the preceding proposition,
Singh's space has the $AP(\infty,0)$-property. Thus we have:

\begin{cor} There is a compact $X\in AR$ with $X\in AP(\infty,0)$ such
that $X\times\I$ is homeomorphic to the Hilbert cube but $X$
contains no proper $ANR$-subspace of dimension $\ge 2$.
\end{cor}

Now, we consider the spaces with piecewise embedding dimension $n$.
According to \cite{km}, a map $h\colon P\to M$ from a finite
polyhedron $P$ is said to be a piecewise embedding if there is a
triangulation $T$ of $P$ such that $h$ embeds each simplex
$\sigma\in T$ and $h(P)$ is an $ANR$. For a space $M$, the piecewise
embedding dimension $ped(M)$ is the maximum $k$ such that for any
$\e>0$ and any map $g\colon P\to M$ from a finite polyhedron $P$
with $\dim P\leq k$ there exists a piecewise embedding $g'\colon
P\to M$ which is $\e$-close to $g$. If $y\in M$, then $ped_y(M)$ is
the maximum of all $ped(U)$, where $U$ is a neighborhood of $y$ in
$M$. Obviously, $ped(M)\leq\min\{ped_y(M):y\in m\}$.

As we noted, by \cite[Propsition 2.1]{km}, every complete
$ANR$-space $M$ with $ped(M)\geq n$ has the $AP(n,0)$-property. But
there are even compact $ANR$'s having the $AP(n,0)$-property with
$ped\leq n-1$. For example, according to Corollary 6.3, any product
$M=\prod_{i=1}^{i=n}M_i$ of dendrites with dense set of endpoints is
an $AP(n,0)$. On the other hand, by \cite[Theorem 3.4]{ys},
$ped(M)\leq n-1$. Next proposition also shows that the property
$ped=n$ is quite restrictive.

\begin{pro}
If $(M,\varrho)$ is a complete space and $ped_y(M)=n+1$, then $y$ is
a $Z_{n}$-point in $M$.
\end{pro}

\begin{proof}
Since $ped_y(M)=n+1$, there exists a neighborhood $U_y$ in $M$ with
$ped(U_y)=n+1$. It is easily seen that if $y$ is a $Z_{n}$-point in
$U_y$, it is also a $Z_{n}$-point in $M$. So, we can suppose that
$U_y=M$, and let $\e>0$ and $g\in C(\I^{n},M)$. We identify $\I^n$
with the set $\{(x_1,..,x_{n+1})\in\I^{n+1}:x_{n+1}=0\}$ and let
$\pi:\I^{n+1}\to\I^n$ be the projection. Since $ped(M)=n+1$, there
exists a triangulation $T$ of $\I^{n+1}$ and a map $g_1\in
C(\I^{n+1},M)$ which is $\e/2$-close to $g\circ\pi$ and $g_1|\sigma$
is an embedding for all $\sigma\in T$. Hence, $g_1^{-1}(y)$ consists
of finitely many points. Moreover, there is $\delta>0$ such that for
any $\delta$-close points $x',x''$ in $\I^{n+1}$ we have
$\varrho(g_1(x'),g_1(x''))<\e/2$. Since each $x\in\I^{n+1}$ is a
$Z_n$-point in $\I^{n+1}$, $g_1^{-1}(y)$ is a $Z_n$-set in
$\I^{n+1}$. Hence, there exists a map $h\colon\I^n\to\I^{n+1}$ such
that $h$ is $\delta$-close to $\displaystyle id_{\I^{n}}$ and
$h(\I^n)\subset\I^{n+1}\backslash g_1^{-1}(y)$. Then $g_1\circ h\in
C(\I^n,M)$ is $\e$-close to $g$ and $y\not\in g_1\circ h(\I^n)$.
\end{proof}

Finally, let us complete the paper with the following question:
\begin{question} Let $X\in AR$ be a compact space such that
$X\in\overline{DD}^{\{2,2\}}\cap AP(\infty,0)$.
Is $X$ homeomorphic to the Hilbert cube?
\end{question}



\end{document}